\documentclass{amsart}

\usepackage{amsmath, amssymb,verbatim,appendix}
\usepackage[mathscr]{eucal}
\usepackage{amscd}
\usepackage{amsthm}
\usepackage{stmaryrd}
\usepackage{enumerate}
\usepackage{comment}
\usepackage[latin1]{inputenc} 
\usepackage{tikz}
\usetikzlibrary{shapes,arrows}
\usepackage{cite}
\usepackage{url}
\usepackage{tikz-cd}
\usepackage[shortlabels]{enumitem}
\usepackage{xfrac}
\usepackage{mathtools}
\usepackage{hyperref}

\makeatletter

\newtheorem{theorem}{Theorem}[section]
\newtheorem{lemma}[theorem]{Lemma}
\newtheorem{corollary}[theorem]{Corollary}
\newtheorem{proposition}[theorem]{Proposition}
\newtheorem{remark}[theorem]{Remark}
\newtheorem{notation}[theorem]{Notation}
\newtheorem{question}{Question}

\newtheorem{observation}[theorem]{Observation}

\newtheorem*{claim}{Claim}
\newtheorem{fact}[theorem]{Fact}

\theoremstyle{definition}

\newtheorem{definition}[theorem]{Definition}

\newtheorem{innercustomthm}{Theorem}
\newenvironment{customthm}[1]
  {\renewcommand\theinnercustomthm{#1}\innercustomthm}
  {\endinnercustomthm}

\newcommand{\abar}{\bar{a}}

\newcommand{\dbar}{\bar{d}}

\newcommand{\calC}{\mathcal{C}}

\newcommand{\calG}{\mathcal{G}}

\newcommand{\calQ}{\mathcal{Q}}

\newcommand{\calU}{\mathcal{U}}

\newcommand{\wt}[1]{\widetilde{#1}}

\newcommand{\id}{\operatorname{id}}
\def\RM{\operatorname{RM}}

\renewcommand{\Im}{\operatorname{Im}}
\newcommand{\ext}{\mathrm{Ext}}

\newcommand{\dl}{\log_{\delta}}
\newcommand{\Aut}{\mathrm{Aut}}
\def\acl{\operatorname{acl}}
\def\dcl{\operatorname{dcl}}
\def\tp{\operatorname{tp}}

\def\hom{\operatorname{Hom}}

\newcommand{\bg}[3]{\Aut_{#1}(#2/#3)}

\def\forkindep{\mathrel{\raise0.2ex\hbox{\ooalign{\hidewidth$\vert$\hidewidth\cr\raise-0.9ex\hbox{$\smile$}}}}}

\title{Splitting differential equations using Galois theory}

\author{Christine Eagles}
\address{Christine Eagles\\
University of Waterloo\\
Department of Pure Mathematics\\
Mathematics \& Computer\\
Waterloo, ON N2L 3G1\\
Canada}
\email{ceagles@uwaterloo.ca}

\author{L\'eo Jimenez}
\address{L\'eo Jimenez\\
The Ohio State University\\
Department of Mathematics\\
Math Tower\\
Columbus, OH 43210-1174\\
United States}
\email{jimenez.301@osu.edu}

\date{\today}


\begin{document}

\begin{abstract}

This article is interested in pullbacks under the logarithmic derivative of algebraic ordinary differential equations. In particular, assuming the solution set of an equation is internal to the constants, we would like to determine when its pullback is itself internal to the constants. To do so, we develop, using model-theoretic Galois theory and differential algebra, a connection between internality of the pullback and the splitting of a short exact sequence of algebraic Galois groups. We then use algebraic group theory to obtain internality and non-internality results.
    
\end{abstract}

\maketitle

\tableofcontents

\section{Introduction}

Differential algebra and model theory have a long history of fruitful interaction, starting with the work of Robinson \cite{robinson1959concept}, where differentially closed fields were introduced. Both disciplines inform each other: the theory of differentially closed fields is a rich source of examples in model theory, and model theoretic techniques can provide useful tools for the study of algebraic differential equations. A reason for this is $\omega$-stability of the theory $\mathrm{DCF}_0$ of differentially closed fields of characteristic zero, which brings all the machinery of finite rank stability theory into the picture.

This article is concerned with one such tool, the \emph{semi-minimal analysis}, and the development of Galois theoretic methods to control it, with an eye towards differential-algebraic applications. Recall that a semi-minimal analysis of a type $p \in S(A)$ is a sequence of $A$-definable maps $p \xrightarrow{f_1} p_1 \xrightarrow{f_2} \cdots \xrightarrow{f_n} p_n$ such that $p_n$, along with each fiber of an $f_i$, is semi-minimal, meaning internal to some rank one type. That this can be done in $\mathrm{DCF}_0$ gives a powerful way to decompose algebraic differential equations.

To make this analysis useful, we need to understand semi-minimal types. For algebraic differential equations, this is done via Zilber's dichotomy for $\mathrm{DCF}_0$: a semi-minimal type is either internal to a locally modular minimal type, and in particular has a rudimentary geometric structure, or is \emph{internal to the constants}. We are interested in the second possibility, of which we give a geometrical account.

A finite dimensional type $p$, over some algebraically closed differential field $k$, is always interdefinable with the generic type of a twisted algebraic vector field $(V,s)$, meaning that $s$ is a morphism over $k$ that is a section of the twisted tangent bundle (or prolongation) $\tau V$. In the autonomous case, i.e. when the derivative on $k$ is trivial, the twisted tangent bundle $\tau V$ is the tangent bundle $TV$. Internality to the constants means that $(V,s)$ is, after possibly taking extensions to a differential field $k < L$, birationally equivalent to the trivial vector field $(\mathbb{A}^{\dim(V)}, 0)$. We also have the close concept of almost internality: the twisted vector field in almost internal to the constants if, after a base change to some $k < K$, there are $\dim(V)$ rational functions on $V$, algebraically independent over $K$, that are constant for the induced derivation. These functions are called \emph{first integrals} in the literature. Note that if we only ask for the existence of one such first integral, we obtain the notion of \emph{non-orthogonality to the constants}. See \cite{moosa2022six} for more on this perspective. 

Even more concretely, the solution set of a differential equation (E) over some differential field $k$ is internal to the constants if there are finitely many solutions $a_1, \cdots , a_n = \abar$ of (E) and a rational function $h$ with coefficients in $k$ such that for any generic solution $a$ of (E), there are constants $c_1, \cdots , c_m$ such that $a = h(\abar, \abar ', \cdots, \abar^{(l)})$ (where the derivative is taken coordinate-wise). This is sometimes called a \emph{superposition principle} in the literature (see \cite{jones1967nonlinear} for example), and generalizes the fact that solution sets of linear differential equations are vectors spaces over the constant field. 

A first step to understand semi-minimal analyses is to study a fibration $f : p \rightarrow f(p)$ with internal fibers, i.e. a definable map with all its fibers internal to a fixed definable set $\calC$. In the case of $\mathrm{DCF}_0$, we will always take $\calC$ to be the field of constants. A starting question would be: suppose that $f(p)$ is also $\calC$-internal, when is $p$ internal to $\calC$? This was one of the motivations of the articles \cite{jaoui2023relative}, \cite{jimenez2019groupoids} and \cite{jin2020internality}. In \cite{jaoui2023relative}, a necessary and sufficient condition is given in the form of \emph{uniform internality}: there are parameters (or a field extension in the concrete case of $\mathrm{DCF}_0$) witnessing internality of all the fibers at once. 

As a guiding example, fix some differential field $k$, some $f \in k(x)$ and consider:
\[ \delta^{m}(y) = f(\delta^{m-1}(y), \cdots, \delta(y),y) \tag{1} \] \label{eqn: base}
\noindent as well as the following system of differential equations:
\[\begin{cases}\tag{2} \label{eqn: pullback}
    \delta^{m}(y) = f(\delta^{m-1}(y), \cdots, \delta(y),y)\\
    \delta(x) = y x
\end{cases}\] 
\noindent On non-zero $x$, there is a logarithmic derivative map $\dl : (x,y) \rightarrow \frac{\delta(x)}{x} = y$ which takes a solution of this system to a solution of (\ref{eqn: base}). It is well-known (and easy to prove) that any fiber of $\dl$ is internal to the constants, and thus this fits into our general setup. We want to understand (\ref{eqn: pullback}) by seeing it as fibered over $\delta^{m}(y) = f(\delta^{m-1}(y), \cdots, \delta(y),y)$. In particular, we want to answer:

\begin{question}\label{qu: internality-log-diff-pullbacks}
    If equation (\ref{eqn: base}) is internal to the constants, is equation (\ref{eqn: pullback}) also (almost) internal to the constants?
\end{question}

From the point of view of vector fields, we know that each fiber of $\dl$ has a first integral, and equation (\ref{eqn: base}) has $m$ first integrals by assumption. What the question is asking is whether from these, we can obtain $m+1$ first integrals to equation (\ref{eqn: pullback}). There is no reason for this to be true in general, and many counterexamples are given in \cite{jin2020internality}. As explained before, this also has direct impact on our ability to recover solutions of equation (\ref{eqn: pullback}) from the data of finitely many fixed solutions.

Note that because arbitrary field extensions are allowed, there is a priori no effective way to check if a given equation is internal to the constants. This is a general phenomenon in model-theoretic differential algebra: many useful concepts involve taking arbitrary field extensions, and controlling which extension one needs to look at has been an important theme. See for example recent work in \cite{freitag2022any} and \cite{freitag2023bounding} on the degree of non-minimality, which led to the development of effective tools for the study of algebraic differential equations. 

Here, we will achieve some effectivity by considering a sufficient condition for internality, introduced in \cite{jin2020internality} and which we will call \emph{splitting}. Going back to the general case of a fibration $f : p \rightarrow f(p)$, we say that $f$ splits if there is an $A$-definable bijection $\iota$ between $p$ and the Morley product $f(p) \otimes r$, for some $\calC$-internal type $r \in S(A)$, and the diagram
\begin{center}
\begin{tikzcd}
    p \arrow[rr, "\iota"] \arrow[dr, "f"]& & f(p) \otimes r \arrow[dl] \\
    & f(p) &
\end{tikzcd}
\end{center}
\noindent commutes. If we replace the bijection $\iota$ by a finite-to-finite correspondence, we obtain \emph{almost splitting}. Note that there is no need to take extra parameters anymore. In the particular case of differential equations, it removes the need to take an arbitrary field extension. One of the pay-offs for understanding splitting is criteria for internality of algebraic differential equations, for example \cite[Theorem B]{jin2020internality}. We will likewise produce such criteria here.

One of the main tools of \cite{jaoui2023relative} and \cite{jin2020internality} is that, under certain conditions, uniform internality and splitting are equivalent. Then, one can use algebraic characterizations of splitting to obtain strong restrictions on when a fibration can be uniformly internal. This is also at the heart of our methods.

Question \ref{qu: internality-log-diff-pullbacks} has been asked for $m = 1$ in \cite{jin2020internality}, but we remove that restriction here. One issue with picking $m > 1$ (and working over possibly non-constant parameters) is that we loose a convenient criteria due to Rosenlicht \cite{rosenlicht1974nonminimality} that tells us exactly when equation (\ref{eqn: base}) is internal to the constants. To the authors' knowledge, no such criteria is known without these assumptions \footnote{Since this article was first written, the authors \cite{eagles2024internality} have obtained a generalization to $m > 1$, over constant parameters.}.

To answer Question \ref{qu: internality-log-diff-pullbacks}, we proceed in two steps, each of independent interest:

\begin{enumerate}[(a)]
    \item develop general model-theoretic tools to determine sufficient conditions for a definable fibration to split, 
    \item give some concrete algebraic characterization of splitting for the logarithmic derivative.
\end{enumerate}

For step (a), we use model-theoretic Galois theory, which makes sense in any $\omega$-stable theory. Given a fixed definable set $\calC$, any $\calC$-internal type comes equipped with a definable group of transformations, called its binding group. In the case of algebraic differential equations, this corresponds to the Galois group associated to Kolchin's strongly normal extensions \cite{kolchin1953galois}, themselves a generalization of the Galois groups of Picard-Vessiot extensions.
As suggested by the anonymous referee, we now state what our methods yield in a general $\omega$-stable context, which is Theorem \ref{theo: int-implies-almost-split-gen} below:

\begin{customthm}{A}\label{theo: A}
    Fix some $\omega$-stable theory $\mathrm{Th}$ eliminating imaginaries, some algebraically closed set of parameters $A$ and some $A$-definable pure algebraically closed field $\calC$. 
    
    Let $p \in S(A)$ be a weakly $\calC$-orthogonal type and $f : p \rightarrow q$ be an $A$-definable map. Assume that:
    \begin{itemize}
        \item $q$ is $\calC$-internal with binding group definably isomorphic to the $\calC$-points of a linear nilpotent group and
        \item for any $a \models p$, the fiber $\tp(a/f(a)A)$ is stationary and $\calC$-internal with binding group isomorphic to the $\calC$-points of a diagonalizable group.
    \end{itemize}
    Then $p$ is almost $\calC$-internal if and only if the map $f$ almost splits.
\end{customthm}

In \cite{jin2020internality}, Jin and Moosa prove a similar theorem in $\mathrm{DCF}_0$, for $q$ generic type of an equation of the form $\delta(x) = f(x)$ with $f \in F(X)$, internal to the constants and with binding group $G_a$, $G_m$ or $G_a \rtimes G_m$. Theorem \ref{theo: A} is therefore a generalization of their result, with the caveat that we now need the binding group to be linear nilpotent (thus excluding the $G_a \rtimes G_m$ case). The result should also hold for linear solvable binding groups, we reserve this question for future work.

To prove Theorem \ref{theo: A}, we make heavy use of binding groups. Consider a fibration $f : p \rightarrow f(p)$, and assume that $p$ is $\calC$-internal. This gives rise to a definable group morphism from the binding group of $p$ to the binding group of $f(p)$. Using the type-definable groupoids of \cite{jimenez2019groupoids}, we connect splitting of the fibration $f$ to definable splitting of the induced morphism between binding groups, thus obtaining necessary and sufficient conditions for splitting of $f$. 

Because we are considering internality to a pure algebraically closed field $\calC$, binding groups over the constants are always definably isomorphic to the constant points of algebraic groups. In particular, we obtain Theorem \ref{theo: A} by using the very constrained structure of nilpotent linear algebraic groups. 

Regarding step (b), concrete characterizations of splitting were given in \cite{jin2020internality} and \cite{jaoui2023relative}, if $k$ is a field of constants. Let $q$ be the generic type of the equation $\delta^{m}(y) = f(\delta^{m-1}(y), \cdots, \delta(y),y)$, splitting was characterized in the following cases:

\begin{itemize}
    \item in \cite[Theorem B]{jin2020internality} for $m = 1$ and $q$ internal to the constants,
    \item in \cite[Theorem 5.2]{jaoui2023relative}, for arbitrary high $m$, but with $q$ orthogonal to the constants. 
\end{itemize}

It follows from \cite[Corollary 5.5]{jaoui2023relative} that the two criteria are in fact identical when $m = 1$. This raises the possibility of a general characterization, free of assumptions on $m$ and $q$. We achieve that goal under a mild additional assumption:

\begin{customthm}{B}\label{theo: B}
    Let $F$ be an algebraically closed differential field, a rational function $f \in F(x_0, \cdots , x_{m-1})$ and $q$ the generic type of $\delta^{m}(y) = f(y, \delta(y), \cdots , \delta^{m-1}(y))$. Let $p := \dl^{-1}(q)$, and suppose that $p(\calU) \not\subset \acl(q(\calU), \calC,F)$. Then $\dl : p \rightarrow q$ is almost split if and only if there are a non-zero $h \in F(x_0, \cdots , x_{m-1})$, some $e \in F$ and some integer $k \neq 0$ such that:
    \[ (kx_0-e)h = \sum\limits_{i=0}^{m-2} \frac{\partial h}{\partial x_i}x_{i+1} + \frac{\partial h}{\partial x_{m-1}}f + \delta^F(h)\]
    \noindent where $\delta^F$ is the unique derivation on $F(x_0, \cdots , x_{m-1})$ such that $\delta^F \vert_F = \delta$ and $\delta^F(x_i) = 0$ for all $i$.
\end{customthm}

This theorem can be considered a simultaneous generalization of the criteria given by \cite[Corollary 5.5]{jaoui2023relative} and \cite[Theorem B]{jin2020internality}, under the extra assumption that $p(\calU) \not\subset \acl(\dl(p)(\calU), \calC,F)$. Let us say a few words about this assumption. The proof of Theorem \ref{theo: B} goes through the auxiliary notion of \emph{product-splitting}: we first show that almost splitting is equivalent to product splitting, and then show that product splitting is equivalent to the second, algebraic property of Theorem \ref{theo: B}. To deduce product-splitting from almost splitting, we need the binding group of the type $r$ witnessing almost splitting to be infinite. The extra assumption is the weakest we could find that implies this. Note that we can prove that the binding group of $r$ is infinite under other extra assumptions, such as when $q$ is orthogonal to the constants (Corollary \ref{cor: prod-split-when-ortho}), or $F$ is a field of constants (Corollary \ref{cor: prod-split-when-constant}).

Together, Theorem \ref{theo: A} and Theorem \ref{theo: B} allow us to prove that some pullbacks under the logarithmic derivative are not internal: we first show they must be split using Theorem \ref{theo: A}, and then use the concrete algebraic characterization of splitting given by Theorem \ref{theo: B} to prove that splitting is impossible. As an example of application, we prove that if $a \models q$ is a generic solution of a linear differential equation, then $\dl^{-1}(q)$ is never $\calC$-internal, provided that the binding group of $q$ is nilpotent and acts transitively on $q$.

Finally, our methods provide examples of internal types $p$ with a non-split fibration $f : p \rightarrow f(p)$. More precisely, using non-split semi-abelian varieties, as well as a bit of algebraic group cohomology, we also answer a question of Jin and Moosa in \cite{jin2020internality} by showing the existence of a $\calC$-internal type $p$ such that the map $\dl : p \rightarrow \dl(p)$ does not split, with the binding group of $\dl(p)$ isomorphic to an elliptic curve:

\begin{customthm}{C}\label{theo: C}
    There exist an algebraically closed field $F < \mathcal{U}$ and a $\calC$-internal type $p \in S_{1}(F)$ such that $\dl : p \rightarrow \dl(p)$ does not almost-split, and the type $\dl(p)$ is strongly minimal and has a binding  group isomorphic to the constant points of an elliptic curve $E$.
\end{customthm}

In particular, Theorem \ref{theo: C} shows that in Theorem \ref{theo: A}, it is necessary to assume that the binding group of the base $q$ is linear.

We now describe the results of the present work in some details. In the preliminary Section \ref{sec: prels}, we first give some results on binding groups. Most of these are well-known, but we single out Lemma \ref{lem: binding iso quotients} as being of particular interest, and, to the authors' knowledge, new. Given any two $\calC$-internal types $p,q \in S(A)$, it shows that their binding groups have a common definable quotient, which encodes the interactions between the two types over the definable set $\calC$. Next, we give a definable version of the Jordan decomposition of a linear algebraic group, mimicking some of the results of \cite{jaoui2022abelian} regarding definable Chevalley and Rosenlicht decompositions. Using this decomposition, we prove that some short exact sequences of definable groups always definably split, which will later yield our automatic splitting result.

Section \ref{section: splitting} consists of general model-theoretic results on splitting. We always consider a definable fibration $f: p \rightarrow f(p)$ with $\calC$-internal fibers, for some fixed definable set $\calC$. In Subsection \ref{subsection: split-prel}, we define splitting of $f$ and connect it to uniform internality, and prove a basic result we will need. In Subsection \ref{subsection: ses}, we introduce the main tool of our work: the connection between splitting of a fibration and splitting of a short exact sequence of definable groups. To do so, we make use of the type-definable groupoid associated to a fibration that was introduced in \cite{jimenez2019groupoids}. In particular, we use the notion of \emph{retractability} of the groupoid introduced in \cite{goodrick2010groupoids}, also reminiscent of descent arguments, see for example \cite[Proposition 4.2]{pillay1997remarks}. If $\calG$ is a groupoid, it is retractable if there is a definable map $r : \mathrm{Ob}(\calG)^2 \rightarrow \mathrm{Mor}(\calG)$ such that for all $a,b,c \in \mathrm{Ob}(\calG)$, we have $r(a,b) \in \mathrm{Mor}(a,b)$ and $r(b,c) \circ r(a,b) = r(a,c)$. We recall and expand on some results of that article, in particular connecting splitting of $f$ and retractability, but under the very restrictive the condition that $f(p)$ is fundamental. The main new result on groupoids is Lemma \ref{lem: ret-when-not-fun}, which gives a sufficient condition to expand this when $f(p)$ is not fundamental: no subgroup of the binding group of $f(p)$ should map definably and with non-trivial image to the binding group of a fiber. We end the section by stating and proving Theorem \ref{theo: A}.

We remark that all the techniques developed in Sections \ref{sec: prels} and \ref{section: splitting} are applicable to any $\omega$-stable theory. The specific example of the theory $\mathrm{CCM}$ of compact complex manifolds comes to mind as another potential area for applications. This would require a good understanding of definable Galois theory in $\mathrm{CCM}$, which has recently been investigated, see \cite{jaoui2023relative} and \cite{jaoui2022abelian} for example. 

In Section \ref{section: logdiff}, we focus on the case of the logarithmic derivative. In subsection \ref{subsection: log-diff-prel}, we prove various results around splitting of the logarithmic derivative. We connect it with another notion of splitting, which we call \emph{product-splitting}, introduced by Jin in his thesis \cite{jin2019logarithmic} and further studied by Jin and Moosa in \cite{jin2020internality}. In \cite[Conjecture 5.4]{jin2019logarithmic}, Jin conjectures that almost splitting and product splitting are equivalent in the specific case where the image of the fibration is internal. We prove that splitting and product-splitting are equivalent under some various additional assumptions, making progress on Jin's conjecture. We then prove an algebraic characterization of product-splitting, which leads us to prove Theorem \ref{theo: B} (Corollary \ref{cor: almost-split-char}).

In Subsection \ref{subsection: nilpotent}, we consider any fibration $\dl: p \rightarrow \dl(p)$, where $p$ is almost internal to the constants, and the binding group of $\dl(p)$ nilpotent. In that case, we obtain automatic splitting by Theorem \ref{theo: A} (Theorem \ref{theo: splits if nilpotent base}). As an example of application, we prove that the logarithmic derivative pullback of the generic type of a linear differential equation with nilpotent binding group is never almost internal to the constants.

Finally, Subsection \ref{subsection: elliptic} is dedicated to proving Theorem \ref{theo: C} (Corollary \ref{cor: jin-moosa-answer}): we construct an example of a non-split fibration by the logarithmic derivative, answering a question of Jin and Moosa in \cite[bottom of page 5]{jin2020internality}. 

\medskip

\textbf{Acknowledgements.} The authors are thankful to Rahim Moosa for multiple conversations on the subject of this article.
We are also grateful to the anonymous referee, whose careful reading of our work lead to numerous improvements, in particular showing the full equivalence between almost internality and almost splitting in Theorem \ref{theo: A}. 

\section{Preliminaries}\label{sec: prels}

\subsection{Internality and binding groups}\label{subsection: binding groups}

In this article, we work, unless specified otherwise, in a sufficiently saturated model $\calU$ of some $\omega$-stable theory $\mathrm{Th}$, eliminating imaginaries. We will make use, without mentioning it, of the fact that in this context, type-definable groups are always definable (see \cite[Corollary 5.19]{poizat2001stable}). Definable will always mean definable over some parameters, and we will always specify the parameters when relevant. When $p$ is a type, we denote $p(\calU)$ the set of realizations of $p$ in $\calU$, and similarly for realizations of definable sets. We will often conflate definable sets and types with their set of realizations in $\calU$. For example, we may say that that a definable group $G$ acts on a type $p$ to mean that it acts on $p(\calU)$. We fix a set of parameters $A$.

We will assume familiarity with geometric stability theory, for which a good reference is \cite{pillay1996geometric}.

Our results will make heavy use of the binding group of an internal type. Most of the literature uses binding groups over a fixed definable set, but we will sometimes need to work over a family of types. We recall the definition of (almost) internality in that case:

\begin{definition}
    Let $\calQ$ be a family of partial types. A stationary type $p \in S(A)$ is (almost) internal to $\calQ$ if there are some $B \supset A$, some $a \models p$, some $c_1, \cdots , c_n$ such that for all $i$ the type $\tp(c_i/B)$ extends a type in $\calQ$, and:

    \begin{itemize}
        \item $a \in \dcl(c_1, \cdots , c_n,B)$ (resp. $\acl$)
        \item $a \forkindep_A B$.
    \end{itemize}
\end{definition}

In particular, if $X$ is any partial type, the type $p$ is (almost) $X$ internal if and only if it is (almost) $\{ X \}$ internal in the sense of the previous definition. 

By a realization of $\calQ$, we mean any tuple realizing some partial type in $\calQ$. We denote the set of realizations of $\calQ$ by $\calQ(\calU)$, or simply $\calQ$ if the context is clear enough. Internality of $p$ is then equivalent to simply stating that there are parameters $B \supset A$ such that $p(\calU) \subset \dcl(\calQ(\calU),B)$. 

In practice, the minimal assumption for the theory to go through is that the family $\calQ$ is $A$-invariant: any automorphism of $\calU$ fixing $A$ fixes the set of realizations of $\calQ$ setwise. This happens, for example, if every partial type in $\calQ$ is over $A$, \emph{which we assume for the rest of this subsection}.

A very useful fact is that if $p$ is $\calQ$-internal, then the parameters $B$ witnessing it can be chosen equal to $A \cup \{ a_1, \cdots , a_n \}$, where $a_1, \cdots , a_n \models p^{(n)}$. We call the $a_i$ a \emph{fundamental system of solutions} of $p$.

Let $\bg{A}{\calU}{\calQ}$ be the group of automorphisms of $\calU$ fixing $A \cup \calQ$ pointwise. A consequence of the previous discussion is that, if $p$ is $\calQ$-internal, the restriction of any $\sigma \in \bg{A}{\calU}{\calQ}$ to $p(\calU)$ is entirely determined by $\sigma(a_1), \cdots , \sigma(a_n)$, where the $a_i$ form a fundamental system. 

In fact, we can consider the group of permutations of $p(\calU)$ coming from the restriction of an element of $\bg{A}{\calU}{\calQ}$, i.e. the set of maps: 
\[ \bg{A}{p}{\calQ} = \left\{ \sigma : p(\calU) \rightarrow p(\calU) \vert \text{ there is } \wt{\sigma} \in \bg{A}{\calU}{\calQ} \text{ with } \wt{\sigma}\vert_{p(\calU)} = \sigma \right\} \text{ .}\]

Fundamental systems of solutions can be used to prove the following classical theorem (see \cite[Theorem 7.4.8]{pillay1996geometric} for a proof):

\begin{fact}\label{fact: binding group}
    Let $\calQ$ be a family of partial types over $A$, and $p \in S(A)$ a complete type, internal to $\calQ$. Then $\bg{A}{p}{\calQ}$ is isomorphic to an $A$-definable group, and its natural action on $p(\calU)$ is relatively $A$-definable.
\end{fact}

We will also denote this definable group $\bg{A}{p}{\calQ}$. Again, if the family $\calQ$ is reduced to a single partial type, this coincides with the usual binding group. 

We will also need the notion of (weak) orthogonality:

\begin{definition}
    The stationary type $p \in S(A)$ is

    \begin{itemize}
        \item weakly orthogonal to $\calQ$ if any realization $a \models p$ is independent, over $A$, from any tuple of realizations of $\calQ$,
        \item orthogonal to $\calQ$ if for any $B \supset A$, any realisation $a \models p\vert_B$ is independent, over $B$, from any tuple of realizations of $\calQ$.
    \end{itemize}
    
\end{definition}

It is well-known that if $p \in S(A)$ is orthogonal to $\calQ$, then it is not almost $\calQ$-internal (unless it is algebraic). On the other hand, a non-algebraic type $p$ can be both $\calQ$-internal and weakly $\calQ$-orthogonal, and we have the well known fact:

\begin{fact}\label{fact: wo implies transit}

A $\calQ$-internal type $p \in S(A)$ is weakly $\calQ$-orthogonal if and only if $\bg{A}{p}{\calQ}$ acts transitively on $p(\calU)$.
    
\end{fact}

For the rest of this subsection, we fix an $A$-definable set $\calC$. We will consider $\calC$-internal types. 

Note that any non-forking extension of a $\calC$-internal type is $\calC$-internal, and we record another well-known fact:

\begin{fact}\label{fact: bindgrp extension}
    Let $p \in S(A)$ be $\calC$-internal. Then for any $B \supset A$,  the binding group $\bg{B}{p\vert_B}{\calC}$ is a definable subgroup of $\bg{A}{p}{\calC}$.
\end{fact}

\begin{proof}
    Pick a Morley sequence $\abar = a_1, \cdots , a_n$ of realizations of $p \vert_B$. It is also a Morley sequence in $p$, and we can pick $n$ large enough so that it is a fundamental system for both $p$ and $p\vert_B$. It is then easy to check that the map $\iota : \bg{B}{p\vert_B}{\calC} \rightarrow \bg{A}{p}{\calC}$ given by sending any $\sigma \in \bg{B}{p\vert_B}{\calC}$ to the (unique) $\iota(\sigma) \in \bg{A}{p}{\calC}$ such that $\iota(\sigma)(\abar) = \sigma(\abar)$ is definable and injective. 
\end{proof}

Recall that a relatively definable map $\sigma$ from a type-definable set $X$ to a type-definable set $Y$ is a relatively definable set $\Gamma(\sigma) \subset X \times Y$ that is the graph of a function from $X$ to $Y$. In the rest of this article, we forgo the adjective relative. If $X,Y$ are type-definable over $A$ and the map $\sigma$ is defined over $A$, then any $\tau \in \Aut_A(\calU)$ gives another map $\tau(\sigma) : X \rightarrow Y$ with its graph being $\tau(\Gamma(\sigma))$. The following easy and well-known observation will be useful:

\begin{observation}[In any theory]\label{obs: aut applied to map}
    Let $X,Y$ be two $A$-type definable sets, and let $\sigma : X \rightarrow Y$ be a definable map. Then for any $\tau \in \Aut_A(\calU)$, we have that $\tau(\sigma) = \tau\vert_Y \circ \sigma \circ \tau^{-1}\vert_X$.
\end{observation}

\begin{proof}
    We have:
    \begin{align*}
        \tau(\Gamma(\sigma)) & = \tau\left(\left\{ (x,\sigma(x)), x \in X \right\}\right) \\
        & = \left\{ (\tau(x),\tau(\sigma(x))), x \in X \right\}
    \end{align*}
    \noindent and thus for any $x \in X$, we see that $\tau(\sigma)(\tau(x)) = \tau(\sigma(x))$. Applying this identity to $\tau^{-1}(x)$, we get:
    \begin{align*}
        \tau(\sigma)(x) & = \tau(\sigma)(\tau(\tau^{-1}(x))) \\
        & = \tau(\sigma(\tau^{-1}(x))) \\
        & = \tau\vert_Y \circ \sigma \circ \tau^{-1}\vert_X (x)
    \end{align*}
\end{proof}

In the rest of this work, we will often forget about which set one needs to restrict $\tau$ to. This will allow us to write $\tau(\sigma) = \tau \circ \sigma \circ \tau^{-1} = \sigma^{\tau^{-1}}$ which is a slight, but harmless, abuse of notation.

We will often apply this result to definable maps coming from binding groups, and we take the opportunity to point out a subtlety. We constantly view $\sigma \in \bg{A}{p}{\calC}$ in two ways: as a definable bijection on $p$, and as an element of the $A$-definable group $\bg{A}{p}{\calC}$. Let us, for now, denote the first as $\sigma\vert_p$ and the second as $\sigma$. What Fact \ref{fact: binding group} tells us is that there is an $A$-definable group action $\mu : \bg{A}{p}{\calC} \times p \rightarrow p$ such that $\mu(\sigma, \cdot)$ is exactly $\sigma\vert_p$ as a bijection on $p$. If $\tau \in \Aut_A(\calU)$ and $a,b \models p$, then:
\begin{align*}
    \mu(\tau(\sigma),\tau(a)) = \tau(b) & \Leftrightarrow \mu(\sigma , a) = b \text{ as } \tau \in \Aut_A(\calU) \\
    & \Leftrightarrow \sigma\vert_p(a) = b \\
    & \Leftrightarrow \tau(\sigma\vert_p)(\tau(a)) = \tau(b) \text{ by the proof of Observation \ref{obs: aut applied to map}}
\end{align*}
In other words, the actions by $\Aut_A(\calU)$ on $\bg{A}{p}{\calC}$, viewed as an $A$-definable group or as a group of definable bijections of $p$, coincide. In the rest of the article, we will identify them and constantly go back-and-forth between the two.

We record the following well-known consequence of Observation \ref{obs: aut applied to map}:

\begin{corollary}\label{cor: auto-normal-binding}
    Let $p \in S(A)$ be a $\calC$-internal type. Any $A$-definable subgroup $H < \bg{A}{p}{\calC}$ is normal. 
\end{corollary}

\begin{proof}
    Let $H$ be such a subgroup, and consider $\tau \in \bg{A}{p}{\calC}$ and $\sigma \in H$. Then $\tau \circ \sigma \circ \tau^{-1} = \wt{\tau}(\sigma)$, where $\wt{\tau}$ is any extension of $\tau$ to $\calU$. As $\wt{\tau}$ fixes $A$ and $H$ is $A$-definable, this implies that $\wt{\tau}(\sigma) \in H$.
\end{proof}

Given two types $p,q \in S(A)$ such that $p$ is $\calC$-internal, we can also form the binding group $\bg{A}{p}{\calC,q}$. Then:

\begin{lemma}\label{lem: bind over family normal}
    The group $\bg{A}{p}{\calC,q}$ is $A$-definably isomorphic to an $A$-definable normal subgroup of $\bg{A}{p}{\calC}$.
\end{lemma}

\begin{proof}
    The group $\bg{A}{p}{\calC,q}$ is $A$-definable, and acts definably on $p$. Viewing $\bg{A}{p}{\calC}$ and $\bg{A}{p}{\calC, q}$ as groups of definable maps acting on $p$, we see that there is an inclusion $\bg{A}{p}{\calC,q} < \bg{A}{p}{\calC}$.
    
    On the definable groups side, fix some Morley sequence $a_1, \cdots , a_n$ in $p$ long enough to be a fundamental system for $p$ over both $\calC$ and $\{ \calC, q \}$. Then we define a map $\iota$ sending any $\sigma \in \bg{A}{p}{\calC , q}$ to the unique element of $\bg{A}{p}{\calC}$ having the same action on the $a_i$ (and thus on $p$). This map is $a_1, \cdots , a_n, A$ definable. 
    
    Consider some $\tau \in \Aut_A(\calU)$ and $\sigma \in \bg{A}{p}{\calC , q}$, it is easy to see, using Observation \ref{obs: aut applied to map}, that $\iota(\tau(\sigma)) = \tau(\iota(\sigma))$. We compute:
    \begin{align*}
        \tau(\iota)(\sigma) & = \tau \circ \iota \circ \tau^{-1}(\sigma) \\
        & = \tau \circ \iota (\tau^{-1}(\sigma)) \\
        & = \tau (\tau^{-1}(\iota(\sigma))) \\
        & = \iota(\sigma)
    \end{align*}
    \noindent so $\tau(\iota) = \iota$. As this holds for any $\tau \in \Aut_A(\calU)$, we obtain that $\iota$ is $A$-definable. So $\bg{A}{p}{\calC , q}$ is identified with an $A$-definable subgroup of $\bg{A}{p}{\calC}$, which is normal by Corollary \ref{cor: auto-normal-binding}.
\end{proof}

We will make use of the following definable version of Goursat's lemma (see \cite[Chapter I, Exercise 5]{lang2012algebra}), the proof of which is left to the reader:

\begin{lemma}\label{lem: goursat}
    Let $G_1$ and $G_2$ be $A$-definable groups, and $H$ an $A$-definable subgroup of $G_1 \times G_2$. Consider the natural projections $\pi_i : G_1 \times G_2 \rightarrow G_i$ for $i = 1,2$, and assume that the maps $\pi_i \vert_H$ are surjective. Let $N_1 = \ker(\pi_2\vert_H)$ and $N_2 = \ker(\pi_1\vert_H)$, which we can identify with $A$-definable normal subgroups of $G_1$ and $G_2$. Then $G_1/N_1$ and $G_2/N_2$ are $A$-definably isomorphic.
\end{lemma}

Here is a consequence regarding the interaction between binding groups:

\begin{lemma}\label{lem: binding iso quotients}
    Let $p,q\in S(A)$ with $p$ and $q$ being $\calC$-internal. Then the quotients $\bg{A}{p}{\calC}/\bg{A}{p}{q,\calC}$ and $\bg{A}{q}{\calC}/\bg{A}{q}{p,\calC}$ are $A$-definably isomorphic.
\end{lemma}

\begin{proof}
    Consider the subgroup $H$ of $\bg{A}{p}{\calC} \times \bg{A}{q}{\calC}$ consisting of pairs $(\sigma_1, \sigma_2)$ having a common extension to an automorphism in $\Aut_A(\calU)$. Let us first show it is $A$-definable. 

    We start by showing it is type-definable over extra parameters. If we fix two fundamental systems of solutions $a_1, \cdots ,a_n$ and $b_1, \cdots , b_m$ of $p$ and $q$, then $(\sigma_1,\sigma_2) \in H$ if and only if $a_1 \cdots a_n b_1 \cdots b_m \equiv_{A\calC} \sigma_1(a_1 \cdots a_n) \sigma_2(b_1 \cdots b_m)$, by stable embeddedness of $\calC$. Still by stable embeddedness, we know that for any tuple $\dbar$, we have $\tp\left(\dbar/ \dcl(\calC A) \cap \dcl(\dbar A) \right) \vdash \tp\left(\dbar/\calC A \right)$. Thus, if we consider the set $\mathfrak{F}$ of $A$-definable functions with domain containing $\tp(a_1, \cdots , a_n, b_1, \cdots , b_m/A)$ and with image in $\dcl(\calC)$, we have that $a_1 \cdots a_n b_1 \cdots b_m \equiv_{A\calC} \sigma_1(a_1 \cdots a_n) \sigma_2(b_1 \cdots b_m)$ if and only if $f(a_1 \cdots a_n b_1 \cdots b_m) = f(\sigma_1(a_1 \cdots a_n) \sigma_2(b_1 \cdots b_m))$ for all $f \in \mathfrak{F}$, which is indeed an $(A,a_1,\cdots , a_n, b_1, \cdots , b_m)$-type definable condition. 
 
    By $\omega$-stability, it must then be an $a_1, \cdots, a_n , b_1 , \cdots , b_m$-definable subgroup. We show that it is $A$-definable by showing it is fixed by any $\tau \in \Aut_A(\calU)$. Let $(\sigma_1, \sigma_2) \in H$. By Observation \ref{obs: aut applied to map}, we know that $\tau(\sigma_i) = \left( \sigma_i \right)^{\tau^{-1}}$. If we let $\wt{\sigma} \in \bg{A}{\calU}{\calC}$ be a common extension of $\sigma_1$ and $\sigma_2$, then the automorphism $\left( \wt{\sigma} \right)^{\tau^{-1}}$ also fixes $\calC \cup A$ pointwise, and extends $(\sigma_1)^{\tau^{-1}}$ and $(\sigma_2)^{\tau^{-1}}$. Thus $(\tau(\sigma_1),\tau(\sigma_2)) \in H$. 

    Let $\pi_1, \pi_2$ be the projection maps from $\bg{A}{p}{\calC} \times \bg{A}{q}{\calC}$ to $\bg{A}{p}{\calC}$ and $\bg{A}{q}{\calC}$, respectively. By Lemma \ref{lem: goursat}, all we need to show is that:

    \begin{enumerate}[(a)]
        \item the restriction of the $\pi_i$ to $H$ are surjective,
        \item $\ker(\pi_2\vert_H) = \bg{A}{p}{\calC, q}$
        \item $\ker(\pi_1 \vert_H) = \bg{A}{q}{\calC, p}$.
    \end{enumerate}

    \noindent where we made the identification of $\ker(\pi_i \vert_H)$ with subgroups of $\bg{A}{p}{\calC}$ and $\bg{A}{q}{\calC}$.

    For (a), we notice that any $\sigma_1 \in \bg{A}{p}{\calC}$ extends to some $\wt{\sigma} \in \bg{A}{\calU}{\calC}$, which restricts to some $\sigma_2 \in \bg{A}{q}{\calC}$, and thus $(\sigma_1, \sigma_2) \in H$, which implies that $\pi_1 \vert_H$ is surjective. We use the same argument for $\pi_2$.

    Statements (b) and (c) are similar, and we only prove (b). Let $(\sigma_1, \id\vert_q)\in\ker(\pi_2\vert_H).$ Then $\sigma_1$ and $\id\vert_q$ have a common extension $\wt{\sigma} \in \bg{A}{\calU}{\calC}$, showing $\sigma_1 \in \bg{A}{p}{\calC, q}$. For the reverse inclusion, let $\sigma_1\in \bg{A}{p}{\calC, q}$, this exactly means that $(\sigma_1, \id \vert_q) \in H$.
\end{proof}

Consider the condition $\bg{A}{p}{\calC} = \bg{A}{p}{\calC, q}$, or equivalently, that any $\sigma \in \bg{A}{p}{\calC}$ extends to $\wt{\sigma} \in \bg{A}{\calU}{\calC}$ fixing $q$ pointwise. An immediate consequence is that this is symmetrical: any element of $\bg{A}{q}{\calC}$ also extends to an automorphism of $\calU$ fixing $p$ pointwise. A sufficient condition for this to happen is for the two binding groups to not have any isomorphic definable quotients:

\begin{corollary}\label{cor: wo from groups}
    If there is no $A$-definable normal subgroups $H \lneq \bg{A}{p}{\calC}$ and $K \lneq \bg{A}{q}{\calC}$ such that $\bg{A}{p}{\calC}/H$ and $\bg{A}{q}{\calC}/K$ are $A$-definably isomorphic, then $\bg{A}{p}{\calC} = \bg{A}{p}{q , \calC}$. In particular, if $p$ is weakly orthogonal to $\calC$, then $p$ is weakly orthogonal to $\{ \calC, q \}$.
\end{corollary}

\begin{proof}
    By Lemma \ref{lem: binding iso quotients}, we obtain that the quotients $\bg{A}{p}{\calC}/\bg{A}{p}{q,\calC}$ and $\bg{A}{q}{\calC}/\bg{A}{q}{p,\calC}$ are $A$-definably isomorphic. By assumption, this implies that $\bg{A}{p}{q,\calC} = \bg{A}{p}{\calC}$. As $p$ is weakly orthogonal to $\calC$, we know that $\bg{A}{p}{\calC}$ acts transitively on $p(\calU)$, and thus so does $\bg{A}{p}{q,\calC}$, which precisely means that $p$ is weakly orthogonal to $\{ q, \calC \}$. 
\end{proof}

Opposed to this condition is isogeny: we say two groups $G_1$ and $G_2$ are (definably) \emph{isogenous} if there are finite subgroups $H_1 < G_1$ and $H_2 < G_2$ such that $G_1/H_1$ and $G_2/H_2$ are (definably) isomorphic. Note that if $G_1$ and $G_2$ are abelian varieties, this coincides with the usual algebraic geometry definition of isogenous. This will be connected to interalgebraicity of types:

\begin{definition}
    Let $p,q \in S(A)$, we say that:
    \begin{itemize}
        \item $p$ and $q$ are \emph{interdefinable} if for any $a \models p$, there is $b \models q$ such that $\dcl(aA) = \dcl(bA)$. Note that by compactness, this is equivalent to the existence of an $A$-definable bijection between $p(\calU)$ and $q(\calU)$.
        \item $p$ and $q$ are \emph{interalgebraic} if for any $a \models p$, there is $b \models q$ such that $\acl(aA) = \acl(bA)$. 
    \end{itemize}
\end{definition}

We obtain the following well-known corollary:

\begin{corollary}\label{cor: interalg implies isogenous}
    Let $p,q \in S(A)$ be two $\calC$-internal types, and suppose that they are interalgebraic over $A$. Then $\bg{A}{p}{\calC}$ and $\bg{A}{q}{\calC}$ are $A$-definably isogenous.
\end{corollary}

\begin{proof}
    By Lemma \ref{lem: binding iso quotients}, it is enough to show that the groups $\bg{A}{p}{\calC, q}$ and $\bg{A}{q}{\calC , p}$ are finite. Pick a fundamental system $a_1, \cdots , a_n$ for $p$. Since $p$ and $q$ are interalgebraic over $A$, there are realizations $b_1, \cdots , b_n$ of $q$ such that $a_i \in \acl(b_i, A)$ for all $i$. Therefore there are only finitely many possibilities for $\sigma(a_i)$ with $\sigma \in \bg{A}{p}{\calC, q}$, for all $i$. Thus there are only finitely many possibilities for the action of $\sigma$ on $p$, implying that $\bg{A}{p}{\calC, q}$ is finite. Similarly we show that $\bg{A}{q}{\calC, p}$ is finite.
\end{proof}

If instead the types are interdefinable, the binding groups are isomorphic.

\subsection{Definable Jordan decomposition}\label{subsection: jordan}

In this subsection, we take inspiration from the definable Chevalley and Rosenlicht decompositions of \cite[Subsection 2.3]{jaoui2022abelian} and produce a definable Jordan decomposition for groups isomorphic to an algebraic group. For this subsection, we pick algebraically closed sets of parameters $A \subset B$, and we assume that $\calC$ is an $A$-definable purely stably embedded algebraically closed field. 

We will extend the Jordan decomposition for linear solvable algebraic groups to groups definably isomorphic to the $\calC$-points of a linear solvable algebraic group. In particular, we prove the existence of a definable unipotent radical, as well as some basic properties. 

Recall that a \emph{torus} is an algebraic group that is isomorphic to $G^n_m$, for some $n$, and a \emph{diagonalizable group} is a group isomorphic to a definable subgroup of $G^n_m$ for some $n$. If a torus is defined over some algebraically closed field $F$, then this isomorphism is defined over $F$ (see \cite[34.3]{humphreys2012linear}). A \emph{unipotent} group is an algebraic group that is isomorphic to an algebraic subgroup of $U_n$, the group of $n \times n$ upper-triangular matrices with $1$'s along the diagonal.

We call a group definably isomorphic to the $\calC$-points of a torus a \emph{definable torus}, and use a similar terminology for definable unipotent, nilpotent, solvable linear algebraic groups and diagonalizable groups.

In this subsection, as well as the rest of this article, we will use, without mentioning it, the fact that for any pure algebraically closed field $\calC$ of characteristic zero, definable groups and algebraic groups coincide, and that any definable map between algebraic groups is a morphism of algebraic groups.

We recall the multiplicative Jordan decomposition, which is a consequence of the Lie-Kolchin Theorem (see \cite[Chapter 19]{humphreys2012linear} for a proof):

\begin{fact}\label{fact: Jordan-decomposition}
    
    Let $G$ be a connected solvable linear algebraic group. There is a split short exact sequence of algebraic groups:
    \[ 1 \rightarrow G_u \rightarrow G \rightarrow T \rightarrow 1\]
    \noindent In particular $G_u \trianglelefteq G$, and moreover, it is the unique maximal unipotent subgroup of $G$ and the group $T$ is a torus. 
    
\end{fact}

We prove that there is a definable Jordan decomposition:

\begin{lemma}\label{lemma: def-jordan-dec}
    Let $G$ be a solvable $B$-definable group, and $f : G \rightarrow \widehat{G}(\calC)$ be a definable isomorphism, where $\widehat{G}$ is a connected linear algebraic group. Then there is a $B$-definable normal subgroup $G_u \trianglelefteq G$ which is the unique maximal subgroup of $G$ definably isomorphic to the $\calC$-points of a unipotent group. We thus obtain a definable exact sequence: 
    \[ 1 \rightarrow G_u \rightarrow G \rightarrow G/G_u \rightarrow 1\]
    \noindent which is definably (maybe over extra parameters) split, and $T := G/G_u$ is definably isomorphic to the $\calC$ points of a torus. 
    
\end{lemma}

The proof is similar to the proof of \cite[Fact 2.8]{jaoui2022abelian}:

\begin{proof}
By Fact \ref{fact: Jordan-decomposition} we have a definable, definably split, short exact sequence $1 \rightarrow \widehat{G}_u \rightarrow \widehat{G} \rightarrow \widehat{T} \rightarrow 1$. Let $G_u=f^{-1}(\widehat{G}_u(\calC))$, then the sequence $1 \rightarrow G_u \rightarrow G \rightarrow G/G_u \rightarrow 1$ is definably isomorphic to the previous sequence of algebraic groups, and thus also definably split. The group $G_u$ is a definable subgroup of $G$ which is definably isomorphic to $\widehat{G}_u(\calC),$  the $\calC$-points of the unipotent group $\widehat{G}_u$.

To show that $G_u$ is unique and maximal, let $H \leq G$ be another definable subgroup, and let $l$ be a definable isomorphism such that $l(H) = \widehat{H}(\calC)$, the $\calC$-points of an unipotent linear algebraic group $\widehat{H}$. Since $\calC$ is purely stably embedded, the map $f \circ l^{-1}$ is a morphism of algebraic groups, and thus $f(H)$ is a unipotent algebraic subgroup of $\widehat{G}(\calC)$. By Fact \ref{fact: Jordan-decomposition}, it must be contained in $\widehat{G}_u(\calC)$, and therefore $H \leq G_u$.

In fact $G_u$ is $B$-definable since any $B$-conjugate of $G_u$ will also be a maximal definable subgroup of $G$ which is definably isomorphic to the $\calC$ points of a unipotent algebraic group. By uniqueness of $G_u$, any such $B$-conjugate is equal to $G_u$, hence $G_u$ is $B$-definable.

Finally, the algebraic group $\widehat{T}$ is a torus, and $T$ is definably isomorphic to $\widehat{T}(\calC)$.
\end{proof}

We call this subgroup the definable unipotent radical of $G$. We show that, just as is the case for linear algebraic groups, it is preserved by definable maps:

\begin{lemma}\label{lem: preserve-def-uni-rad}
    Let $\pi: G \rightarrow K$ be a map between solvable $B$-definable groups and let $f : G \rightarrow \widehat{G}(\calC)$ and $h: K \rightarrow \widehat{K}(\calC)$ be definable isomorphisms, where $\widehat{G}$ and $\widehat{K}$ are connected linear algebraic groups. Then $\pi(G_u) < K_u$
\end{lemma}

\begin{proof}
    We obtain, by stable embeddedness of $\calC$, a $\calC \cap B$-definable morphism of linear algebraic groups $\widehat{\pi} = h \circ \pi \circ f^{-1} : \widehat{G}(\calC) \rightarrow \widehat{K}(\calC)$, and we must have, by \cite[Section 15.3]{humphreys2012linear} that $\widehat{\pi}(\widehat{G}_u) < \widehat{K}_u$. Applying $h^{-1}$, we obtain $\pi \circ f^{-1}(\widehat{G}_u) = h^{-1}(\widehat{K}_u)$. By construction of the definable unipotent radical, this gives us $\pi(G_u) < K_u$.
\end{proof}

If the group is nilpotent, we obtain more:

\begin{lemma}\label{lem: def-prod-nilp}
    Let $G$ be a nilpotent $B$-definable group, and $f : G \rightarrow \widehat{G}(\calC)$ be a definable isomorphism, where $\widehat{G}$ is a connected linear algebraic group. Then there are $B$-definable normal subgroups $G_u \trianglelefteq G$ and $T \trianglelefteq G$ which are the unique maximal subgroup of $G$ definably isomorphic to the $\calC$-points of a unipotent group, respectively a torus. Moreover $G$ is $B$-definably isomorphic to $G_u \times T$.
\end{lemma}
\begin{proof}
    Here $\widehat{G}$ is a nilpotent connected linear algebraic group, and thus by \cite[Proposition 19.2]{humphreys2012linear} splits into a product $\widehat{G_u} \times \widehat{T}$ for some unipotent $\widehat{G_u}$ and torus $\widehat{T}$. These groups are the maximal unipotent subgroup and torus of $\widehat{G}$. Using this, we can prove, as was done in Lemma \ref{lemma: def-jordan-dec}, that their preimages under $f$ are normal $B$-definable subgroups of $G$, and it then follows that $G$ is $B$-definably isomorphic to $G_u \times T$.
\end{proof}
In our application, we will obtain a short exact sequence that goes the other way: its kernel is a torus, and its image nilpotent. Having both of these forces splitting as a direct product:

\begin{lemma}\label{lem: nilp-splitting}
    Let $G$ be a $B$-definable connected group, and $f : G \rightarrow \widehat{G}(\calC)$ a definable isomorphism, where $\widehat{G}$ is an algebraic group. Suppose that there is a short exact sequence:
    \[ 1 \rightarrow K \rightarrow G \xrightarrow{\pi} G/K \rightarrow 1\]
    \noindent where $K$ is a definable diagonalizable group, and $G/K$ is definably nilpotent linear. Then $G$ is nilpotent and $\widehat{G}$ is linear. If moreover $K$ is a torus and the maximal tori of $G$ and $G/K$ are $B$-definably isomorphic to the $\calC$-points of tori, then the short exact sequence is $B$-definably split and $G = K \times G/K$.
\end{lemma}

\begin{proof}
    As an extension of a solvable group by an abelian group, the group $G$ is solvable, and $\widehat{G}$ is linear as an extension of a linear group by a linear group. By \cite[Proposition 19.4]{humphreys2012linear}, as $K$ is a definable diagonalizable group, this implies that $C_G(K) = N_G(K)$, and as $K$ is normal, is must be central. Therefore $G$ is nilpotent because it is a central extension of a nilpotent group. 

    Now suppose that $K$ is a torus. By Lemma \ref{lem: def-prod-nilp}, we know that $G$ is $B$-definably isomorphic to the product of its maximal unipotent subgroup and torus, and the same is true for $G/K$. Consider the following diagram, where all vertical and horizontal short sequences are exact and the vertical sequences are given by Lemma \ref{lem: def-prod-nilp}:
    \begin{center}
    \begin{tikzcd}
    & & 1 \arrow[d] & 1 \arrow[d] &  \\
    & & G_u \arrow[d, "\iota"] & (G/K)_u \arrow[d] &  \\
    1 \arrow[r] & K \arrow[r] & G \arrow[r, "\pi"] \arrow[d] & G/K \arrow[r] \arrow[d]& 1  \\
    & & T \arrow[d] & S \arrow[d]&  \\
    & & 1 & 1 &  \\
    \end{tikzcd}
    \end{center}

\noindent The map $\pi \circ \iota$ has image a subgroup of $(G/K)_u$. Since $K$ is a definable torus and $G_u$ is definably unipotent, we can show that $K \cap G_u = \{ \id \}$, thus $\pi \circ \iota$ is injective. Note that its image must be contained in $(G/K)_u$ by Lemma \ref{lem: preserve-def-uni-rad}. Moreover, as $\pi$ is surjective, the Morley rank of $(G/K)_u$ (which equals the dimension of the linear algebraic group it is isomorphic to) must be smaller or equal to the Morley rank of of $G_u$. Thus the map $\pi \circ \iota$ must be bijective. 

From the diagram, we see that there is an induced short exact sequence:
\[1 \rightarrow K \rightarrow T \rightarrow S \rightarrow 1\]
\noindent By assumption, the maximal tori $T$ and $S$ are $B$-definably isomorphic to the $\calC$-points of tori $\widehat{T}$ and $\widehat{S}$, which are by stable embeddedness defined over $\calC \cap B$. Using the machinery of characters of $d$-groups (see \cite[Section 16.2]{humphreys2012linear}), we can prove that the induced map $\widehat{T} \rightarrow \widehat{S}$ is split. Moreover, because $B \cap \calC$ is algebraically closed, by \cite[Proposition 3.2.12]{springer1998linear}, it is $B \cap \calC$-definably split. Therefore the induced short exact sequence from $T$ to $S$ is $B$-definably split, implying that $T$ is $B$-definably isomorphic to $K \times S$.

Thus we have obtained that $G$ and $G/K$ are $B$-definably isomorphic to $G_u \times S \times K$ are $G_u \times S$, respectively. We can conclude immediately from this. 
\end{proof}

We conclude this section with a lemma on decomposition of actions of definably nilpotent linear algebraic groups:

\begin{lemma}\label{lem: nilp-action-dec}
    Let $G = G_u \times T$ be a definable nilpotent linear algebraic group, where $G_u$ is its unipotent radical and $T$ its maximal torus. Consider a definable group action of $G$ on some definable set $X$. Then for all $a \in X$, we have $G_u a \cap T a = \{ a \}$.
\end{lemma}

\begin{proof}
Let $\pi_1$ and $\pi_2$ be the projections on $G_u$ and $T$. We first prove:

    \begin{claim}
        If $H < G_u \times T$ is a definable subgroup, then $H = \pi_1(H) \times \pi_2(H)$.
    \end{claim}

    \begin{proof}
        By Lemma \ref{lem: goursat} applied to the group $\pi_1(H) \times \pi_2(H)$, there is a definable isomorphism between $\pi_1(H)/N_1$ and $\pi_2(H)/N_2$, where $N_1 = \ker(\pi_2 \vert_H)$ and $N_2 = \ker(\pi_1 \vert_H)$ are normal definable subgroups of $H$, which we can identify with subgroups of $\pi_1(H)$ and $\pi_2(H)$. But the group $\pi_1(H)/N_1$ is definably unipotent and the group $\pi_2(H)/N_2$ is a quotient of a definable torus, so there cannot be any definable isomorphism between them, unless they are both trivial. Therefore $\pi_1(H) = N_1$ and $\pi_2(H) = N_2$, which is equivalent to $H = \pi_1(H) \times \pi_2(H)$.
    \end{proof}
    
    The action of $G$ on $X$ is definably isomorphic to its action on $G/H$, where $H < G$ is a definable subgroup. By the claim, we know that $H = \pi_1(H) \times \pi_2(H)$, and thus the action is definably isomorphic to the action of $G =  G_u \times T$ on $G_u/\pi_1(H) \times T/\pi_2(H)$. In particular for any $a \in X$, we have $G_ua \cap Ta = \{ a \}$.
\end{proof}

\begin{remark}
    We have not used that the ambient theory is $\omega$-stable in this subsection, and the reader can check that our results go through in any theory eliminating imaginaries, as long as $\calC$ is a purely stably embedded algebraically closed field.
\end{remark}

\section{Splitting}\label{section: splitting}

\subsection{Definable fibrations and uniform internality}\label{subsection: split-prel}

In this section, we are still working in an $\omega$-stable theory $T$ eliminating imaginaries. We fix some algebraically closed set of parameters $A$, and some $A$-definable set $\calC$. 

If $p \in S(A)$, by an $A$-definable map $f : p \rightarrow f(p)$, what we mean is an $A$-definable map $f$ with domain containing $p$, which restricts into a relatively definable map on $p(\calU)$. The image $f(p(\calU))$ is then the set of realisations of a type over $A$, which we denote $f(p)$. We will mainly be concerned with definable \emph{fibrations}:

\begin{definition}
    Let $p \in S(A)$. An $A$-definable map $f : p \rightarrow f(p)$ is a fibration if for any $a \models p$, the type $\tp(a/f(a)A)$ is stationary.
\end{definition}

Note that the set of realizations of $\tp(a/f(a)A)$ is exactly the fiber $f^{-1}(\{ f(a) \}) \cap p(\calU)$. We introduce a notation for these types:

\begin{notation}
    If $p \in S(A)$ and $f : p \rightarrow f(p)$ is an $A$-definable map, then for any $f(a) \models f(p)$, we denote $p_{f(a)} : = \tp(a/f(a)A)$.
\end{notation}

More specifically, we will mostly be interested in fibrations such that each fiber $p_{f(a)}$ is almost internal to some fixed definable set $\calC$. These are called \emph{relatively internal} fibrations in \cite{jaoui2023relative}. We are interested in $\calC$-internality of $p$, under the assumption that $f(p)$ is $\calC$-internal. An answer is given by \emph{uniform internality}, introduced in \cite{jaoui2023relative}:

\begin{definition}
    Let $p \in S(A)$, and $f : p \rightarrow f(p)$ a fibration. The fibration $f : p \rightarrow f(p)$ is said to be uniformly $\calC$-internal (resp. uniformly almost $\calC$-internal) if there is a tuple $e$ such that for some (any) $a \models p$, independent from $e$ over $A$, we have that $a \in \dcl(e,f(a),\calC)$ (resp. $\acl$).
\end{definition}

\noindent This implies that the fibers $p_{f(a)}$ are (almost) $\calC$-internal, but is strictly stronger. The following is easy to prove (see \cite[Proposition 3.16]{jaoui2022abelian}):

\begin{fact}
     Let $p \in S(A)$, and $f : p \rightarrow f(p)$ a fibration with almost $\calC$-internal fibers. The type $p$ is almost $\calC$-internal if and only if $f(p)$ is almost $\calC$-internal and $f$ is uniformly almost $\calC$-internal.
\end{fact}

In this article, we are interested in a strictly stronger condition implying internality, which we call \emph{splitting} (and is called triviality in \cite{jaoui2023relative}):

\begin{definition}
    Let $f : p \rightarrow f(p)$ be an $A$-definable fibration with almost $\calC$-internal fibers. We say that $f$ splits (resp. almost splits) if there is an almost $\calC$-internal type $r \in S(A)$ such that for any $a \models p$, there is $b \models r$ such that $(f(a),b)$ and $a$ are interdefinable (resp. interalgebraic) over $A$, and $f(a) \forkindep_A b$. 
\end{definition}

Graphically, splitting means that there is an $A$-definable bijection between $p$ and $f(p) \otimes r$ such that the diagram:

\begin{center}
\begin{tikzcd}
    p \arrow[rr] \arrow[dr, "f"]& & f(p) \otimes r \arrow[dl] \\
    & f(p) &
\end{tikzcd}
\end{center}

\noindent commutes. It is key to our methods to be able to replace almost internal types with internal ones. To do so, the following is useful (see \cite[Lemma 3.6]{jin2020internality} for a proof):
\begin{fact}\label{fact: almost int to int}
    If $p = \tp(a/A)$ is an almost $\calC$-internal type, then there is some $b \in \dcl(aA)$ such that $a \in \acl(bA)$ and $\tp(b/A)$ is $\calC$-internal. In particular, any almost $\calC$-internal type is interalgebraic with a $\calC$-internal type. 
\end{fact}

Applying Fact \ref{fact: almost int to int}, it is easy to prove:

\begin{proposition}\label{pro: reduce to int}
    Let $f : p \rightarrow f(p)$ be an almost split $A$-definable fibration, and let $r$ be the type witnessing it. Then there is a $\calC$-internal type $\widetilde{r} \in S(A)$ such that for any $ a \models p$, there is $b \models \wt{r}$ such that $a$ and $(b,f(a))$ are interalgebraic over $A$.
\end{proposition}

\subsection{Groupoids, splitting and short exact sequences}\label{subsection: ses}

Recall that definable maps give rise to definable morphisms of binding groups (see for example \cite[Lemma 3.1]{jin2020internality}). More precisely: 

\begin{fact}\label{fact: f-tilde}

If $p \in S(A)$ is $\calC$-internal and $f : p \rightarrow f(p)$ is an $A$-definable map, then $f(p)$ is also $\calC$-internal and there is a surjective $A$-definable group homomorphism $\widetilde{f} : \Aut_A(p/\calC) \rightarrow \Aut_A(f(p)/\calC)$ such that for any $\sigma \in \bg{A}{p}{\calC}$ and $a\models p$, we have:
\[ \wt{f}(\sigma)(f(a)) = f(\sigma(a)) \text{ .}\]
\end{fact}

Our main tool will be a connection between the splitting of a fibration and the definable splitting of the definable short exact sequence arising from it. 

We will use the type-definable groupoid associated to the fibration $f$, introduced by the second author in \cite{jimenez2019groupoids} and \cite{jimenez2020internality}. 
We briefly recall the relevant facts and definitions and direct the reader to the aforementioned works for more details.

Let $f : p \rightarrow f(p)$ be an $A$-definable fibration with $\calC$-internal fibers. For any $f(a), f(b) \models f(p)$, we define the set of morphisms $\mathrm{Mor}(f(a),f(b))$ to be the set of bijections from $p_{f(a)}$ to $p_{f(b)}$ which extend to an automorphism of $\calU$ fixing $A \cup \calC$ pointwise. We define the groupoid associated to $f$, denoted $\mathrm{Gd}_A(p,f/\calC)$, to be the groupoid with set of objects $f(p)$ and morphisms the sets $\mathrm{Mor}(f(a),f(b))$. This groupoid acts naturally on $p$ in the following way: by construction, any element of $\mathrm{Mor}(f(a),f(b))$ is a bijection from $p_{f(a)}$ to $p_{f(b)}$, and we define its (partial) action on $p$ to be that bijection. The following is \cite[Theorem 1.3]{jimenez2019groupoids}:

\begin{theorem}
    The groupoid $\mathrm{Gd}_A(p,f/\calC)$ is isomorphic to an $A$-type-definable groupoid, and its natural action on $p$ is relatively $A$-definable.
\end{theorem}

We recall the definition of retractability from \cite{goodrick2010groupoids}:

\begin{definition}
    The groupoid $\mathrm{Gd}_A(p,f/\calC)$ is \emph{retractable} if there is an $A$-definable map $r: f(p)^2 \rightarrow \mathrm{Gd}_A(p,f/\calC)$ such that for all $a,b,c \models f(p)$:
    \begin{itemize}
        \item $r(f(a),f(b)) \in \mathrm{Mor}(f(a),f(b))$,
        \item $r(f(b),f(c)) \circ r(f(a),f(b)) = r(f(a),f(c))$,
    \end{itemize}
\end{definition}

Roughly speaking, retractability means that we can pick, uniformly definably, morphisms between the fibers of $f$, in a way that is compatible with composition. Note that it implies, in particular, that $\mathrm{Gd}_A(p,f/\calC)$ is \emph{connected}: any two of its objects have a morphism between them. By Fact \ref{fact: wo implies transit} this is equivalent to $f(p)$ being weakly $\calC$-orthogonal. 

We will use the following:

\begin{lemma}\label{lem: split-iff-retract}
    Let $f: p \rightarrow f(p)$ be an $A$-definable fibration with $\calC$-internal fibers. The groupoid $\mathrm{Gd}_A(p,f/\calC)$ is retractable if and only if there exists an $A$-definable map $\pi : p \rightarrow \pi(p)$ such that:
    \begin{itemize}
        \item $f(p)$ is weakly orthogonal to $\{ \calC, \pi(p) \}$,
        \item $\pi \times f : p \rightarrow \pi(p) \times f(p)$ is a definable bijection.
    \end{itemize}
    In particular retractability of $\mathrm{Gd}_A(p,f/\calC)$ implies splitting of $f$.
\end{lemma}

For a proof, going from retractability to splitting is given by\cite[Proposition 4.3]{jimenez2019groupoids}. However as pointed out by the anonymous referee, the weak orthogonality of $f(p)$ to $\{ \calC , \pi(p)\}$ is not proved there. The other direction is given in the proof of \cite[Theorem 3.3.5]{jimenez2020internality}, and is similar to the right to left direction of the next lemma. For these reasons, we do not give a full proof, but only a brief sketch of the left to right implication below.

Assuming $\mathrm{Gd}_A(p,f/ \calC)$ is retractable, we have a family of map $r(f(a),f(b))$. We define an equivalence relation $E$ on $p$ as $E(x,y)$ if and only if $r(f(x),f(y))(x) = y$. The compatibility conditions of retractability exactly mean that this is an equivalence relation. The map $\pi$ is the quotient map of $E$. We need to show that $f(p)$ is weakly orthogonal to $\{ \calC, \pi(p) \}$, or equivalently, that $\bg{A}{\calU}{\calC, \pi(p)}$ acts transitively on $f(p)$. Pick any $f(a), f(b) \models f(p)$ and let $\sigma$ be any extension of $r(f(a),f(b))$ to an element of $\bg{A}{\calU}{\calC}$. 

Since $E$ is $A$-definable, we know that $\sigma$ permutes the classes of $E$. Any $E$-class has a representant $a' \models p_{f(a)}$, and we have that $E(a', \sigma(a'))$. Therefore $\sigma$ fixes the class of $a'$ pointwise. As this is true for all classes, we get that $\sigma$ fixes $\pi(p)$ pointwise, and therefore restricts to an element of $\bg{A}{f(p)}{\calC, \pi(p)}$. This shows that $f(p)$ is weakly orthogonal to $\{ \calC, \pi(p) \}$, and the second condition is easy to show (and proved in \cite[Proposition 4.3]{jimenez2019groupoids}).

When $f(p)$ is fundamental, the connection between splitting of the short exact sequence and retractability is straightforward. This is proven in \cite[Theorem 4.10]{jimenez2019groupoids}, but we include a considerably streamlined proof:

\begin{lemma}\label{lem: fund-direct-split-retract}
    Let $p \in S(A)$ be a $\calC$-internal type and $f : p \rightarrow f(p)$ be an $A$-definable fibration. Suppose $f(p)$ is weakly $\calC$-orthogonal and fundamental. Then the short exact sequence 
    \[ 1 \rightarrow K \rightarrow \Aut_A(p/\calC) \xrightarrow{\wt{f}} \Aut_A(f(p)/\calC) \rightarrow 1\]
    \noindent is $A$-definably split if and only if the groupoid $\mathrm{Gd}_A(p,f/\calC)$ is retractable (which implies $f$ splits). 
    
    Moreover in that case, the binding group $\bg{A}{p}{\calC}$ is $A$-definably isomorphic to $\bg{A}{f(p)}{\calC} \times \bg{A}{\pi(p)}{\calC}$, where $\pi$ is the map of Lemma \ref{lem: split-iff-retract}. 
\end{lemma}

\begin{proof}
    Suppose first that the short exact sequence is $A$-definably split. Then there exists an $A$-definable map $s : \bg{A}{f(p)}{\calC} \rightarrow \bg{A}{p}{\calC}$ such that $\wt{f} \circ s = \id$. For any $f(a),f(b) \models f(p)$, because $f(p)$ is fundamental, there exists a unique $\sigma_{f(a),f(b)} \in \bg{A}{f(p)}{\calC}$ such that $\sigma_{f(a),f(b)}(f(a)) = f(b)$. It is easy to see that the family of maps $\{ s(\sigma_{f(a),f(b)})\vert_{p_{f(a)}} : f(a),f(b) \models f(p)\}$ witnesses retractability of $\mathrm{Gd}_A(p,f/\calC)$.
    
    Conversely, suppose that $\mathrm{Gd}_A(p,f/\calC)$ is retractable, and consider the map $\pi$ given by Lemma \ref{lem: split-iff-retract}. The maps $f$ and $\pi$ induce an $A$-definable group morphism $\wt{\pi} \times \wt{f} : \bg{A}{p}{\calC} \rightarrow \bg{A}{\pi(p)}{\calC} \times \bg{A}{f(p)}{\calC}$. As $p$ is interdefinable with $\pi(p) \times f(p)$, this map is injective. 

    As $f(p)$ is fundamental, the action of $\bg{A}{f(p)}{\calC}$ is free and transitive. As $f(p)$ is weakly orthogonal to $\left\{ \pi(p), \calC \right\}$, so is the action of $\bg{A}{f(p)}{\calC, \pi(p)}$, and therefore we must have $\bg{A}{f(p)}{\calC, \pi(p)} = \bg{A}{f(p)}{\calC}$, and by Lemma \ref{lem: binding iso quotients} we also obtain $\bg{A}{\pi(p)}{\calC, f(p)} = \bg{A}{\pi(p)}{\calC}$.

    It is easy to see that $\bg{A}{\pi(p)}{\calC, f(p)} \times \bg{A}{f(p)}{\calC, \pi(p)} < \Im(\wt{\pi} \times \wt{f})$: given $(\sigma, \tau)$ in the former group, we can extend each to $\widehat{\sigma}, \widehat{\tau} \in \bg{A}{p}{\calC}$, and then $(\wt{\pi} \times \wt{f})(\widehat{\sigma} \widehat{\tau}) = (\sigma, \tau)$. Therefore the map $\wt{\pi} \times \wt{f}$ is an $A$-definable isomorphism between $\bg{A}{p}{\calC}$ and $\bg{A}{\pi(p)}{\calC} \times \bg{A}{f(p)}{\calC}$.
\end{proof}

Note that even when $f(p)$ is not fundamental, if there is an $A$-definable section $s$, by Corollary \ref{cor: auto-normal-binding}, the group $\Im(s)$ is an $A$-definable normal subgroup of $\bg{A}{p}{\calC}$, and this group is $A$-definably isomorphic to $\bg{A}{f(p)}{\calC} \times \Im(s)$. 

There is one important case where the type $f(p)$ is fundamental: if its binding group is abelian. Indeed, the action of $\bg{A}{f(p)}{\calC}$ is faithful, and all stabilisers are conjugate, hence equal. So if some $\sigma \in \bg{A}{f(p)}{\calC}$ fixes some $f(a) \models f(p)$, then its fixes every realisation of $f(p)$, and thus must be the identity by faithfullness. We will use this to deal with binding groups isomorphic to tori.

The assumption that $f(p)$ is fundamental can always be obtained, without changing the short exact sequence, by taking a high enough Morley power of $p$ (see \cite[Fact 2.4]{jaoui2022abelian}). However, applying Lemma \ref{lem: fund-direct-split-retract} would yield a splitting of a Morley power of $p$, which is not what we want. Note that the right to left direction of the proof only uses that $f(p)$ is fundamental to obtain the equality $\bg{A}{f(p)}{\calC, \pi(p)} = \bg{A}{f(p)}{\calC}$. We could therefore replace the assumption that $f(p)$ is fundamental by this and obtain the same implication.

Outside of the abelian case, we cannot expect the type $f(p)$ to be fundamental. In particular, we will need something more to deal with unipotent binding groups. The new results of this section give conditions allowing us to bypass this restriction. 

\begin{definition}
    Let $p \in S(A)$ be a $\calC$-internal type and $f : p \rightarrow f(p)$ be an $A$-definable fibration. Suppose $f(p)$ is weakly $\calC$-orthogonal and that the short exact sequence 
    \[1 \rightarrow K \rightarrow \Aut_A(p/\calC) \xrightarrow{\wt{f}} \Aut_A(f(p)/\calC) \rightarrow 1\]
    \noindent is $A$-definably split via some $A$-definable section $s$.

    Let $H$ be a definable subgroup of $\bg{A}{f(p)}{\calC}$. We say that $H$ satisfies $(\bigstar)$ if for any $f(a) \models f(p)$ and $\sigma \in H$, if $\sigma(f(a)) = f(a)$, then $s(\sigma) \vert_{p_{f(a)}} = \id$. 
\end{definition}

Let us point out what $(\bigstar)$ means. If $\sigma$ fixes $f(a)$, then $s(\sigma)$ fixes $p_{f(a)}$, the fiber above $f(a)$, as a set. This assumption says that it in fact fixes $p_{f(a)}$ pointwise. We obtain:

\begin{lemma}\label{lem: ret-when-not-fun}
    Let $p \in S(A)$ be a $\calC$-internal type and $f : p \rightarrow f(p)$ be an $A$-definable fibration. Suppose $f(p)$ is weakly $\calC$-orthogonal and that the short exact sequence $1 \rightarrow K \rightarrow \Aut_A(p/\calC) \xrightarrow{\wt{f}} \Aut_A(f(p)/\calC) \rightarrow 1$ is $A$-definably split via some $A$-definable section $s$. If $\bg{A}{f(p)}{\calC}$ satisfies $(\bigstar)$, then $\mathrm{Gd}_A(p,f/\calC)$ is retractable (and in particular $f$ splits).
\end{lemma}

\begin{proof}
    Let $f(a),f(b) \models f(p)$ and $\sigma, \tau \in \bg{A}{f(p)}{\calC}$ such that $\sigma(f(a)) = f(b) = \tau(f(a))$. Then $s(\sigma)\vert_{p_{f(a)}}$ and $s(\tau)\vert_{p_{f(a)}}$ are in $\mathrm{Mor}(f(a),f(b))$. By assumption $(s(\sigma)\vert_{p_{f(a)}})^{-1} \circ s(\tau)\vert_{p_{f(a)}} = s(\sigma \circ \tau^{-1})\vert_{p_{f(a)}} = \id \vert_{p_{f(a)}}$, so $s(\sigma)\vert_{p_{f(a)}} = s(\tau)\vert_{p_{f(a)}}$. Therefore, we can define, for any $f(a),f(b) \models f(p)$, the map $r(f(a),f(b)) = s(\sigma)_{p_{f(a)}}$, for any $\sigma \in \bg{A}{f(p)}{\calC}$ with $\sigma(f(a)) = f(b)$. This is well-defined by the previous discussion, and it is straightforward to show that it witnesses retractability of $\mathrm{Gd}_A(p,f/\calC)$.
\end{proof}

Note that if $f(p)$ is fundamental, then any $\sigma \in \bg{A}{f(p)}{\calC}$ fixing some $f(a) \models f(p)$ must be the identity. Therefore in that case, we always have (\emph{$\bigstar$}). If $f(p)$ is not fundamental, we find the following sufficient condition for (\emph{$\bigstar$}):

\begin{corollary}\label{cor: no-map-implies-star}
    Under the same assumptions as Lemma \ref{lem: ret-when-not-fun}, if for some (any) $f(a) \models f(p)$, there is no non-trivial (i.e. with image not equal to the trivial subgroup) definable morphism from any definable subgroup of $\bg{A}{f(p)}{\calC}$ to $\bg{Af(a)}{p_{f(a)}}{\calC}$ then assumption \emph{($\bigstar$)} holds. In particular $f$ splits.
\end{corollary}

\begin{proof}
    Fix $f(a) \models f(p)$, let $H$ be the definable subgroup of elements $\bg{A}{f(p)}{\calC}$ fixing $f(a)$. Then $s$ gives rise to a definable map $s' : H \rightarrow \bg{f(a)A}{p_{f(a)}}{\calC}$, which must have trivial image. This is exactly saying that if $\sigma(f(a)) = f(a)$, then $s(\sigma)$ fixes $p_{f(a)}$ pointwise. 
\end{proof}

To go from unipotent and tori to nilpotent, we need to remark that (\emph{$\bigstar$}) is preserved under direct product. More precisely:

\begin{lemma}\label{lem: starprod}
    Let $p \in S(A)$ be a $\calC$-internal type and $f : p \rightarrow f(p)$ be an $A$-definable fibration. Suppose $f(p)$ is weakly $\calC$-orthogonal and that the short exact sequence $1 \rightarrow K \rightarrow \Aut_A(p/\calC) \xrightarrow{\wt{f}} \Aut_A(f(p)/\calC) \rightarrow 1$ is $A$-definably split via some $A$-definable section $s$. Suppose that $\bg{A}{f(p)}{\calC}$ is $A$-definably isomorphic to a direct product $H_1 \times H_2$, that for all $f(a) \models f(p)$ we have $H_1 f(a) \cap H_2 f(a) = \{ f(a) \}$, and that the induced actions of $H_1$ and $H_2$ on $f(p)$ both satisfy \emph{($\bigstar$)}. Then the action of $\bg{A}{f(p)}{\calC}$ on $f(p)$ satisfies \emph{($\bigstar$)}. 
\end{lemma}

\begin{proof}
    We can write any $\sigma \in \bg{A}{f(p)}{\calC}$ as $\sigma_1 \sigma_2$ with $\sigma_i \in H_i$. Suppose that $\sigma(f(a)) = f(a)$ for some $f(a) \models f(p)$. Then $\sigma_1(f(a)) = \sigma_2^{-1}(f(a))$, therefore $\sigma_1(f(a)) = \sigma_2(f(a)) = f(a)$. By our (\emph{$\bigstar$}) assumptions, we obtain that $s(\sigma_i) \vert_{p_{f(a)}} = \id$, and therefore $s(\sigma)\vert_{p_{f(a)}} = \id$.
\end{proof}

Let us say a word about the assumptions of these lemmas before continuing. First, the assumption that the short exact sequence is $A$-definably split cannot be weakened: a counterexample is given in Jin's thesis \cite[Example 5.29]{jin2019logarithmic}. Here is a quick presentation of it: fix $\calU \models \mathrm{DCF}_0$ a sufficiently saturated model, and $\calC$ its field of constants. If $s$ is a differential transcendental, Jin considers the $\calC$-internal types $p$ generic of $\delta^{-1}(s)$ and $q = \delta^{-1}(p)$. He shows that the short exact sequence obtained is isomorphic to an extension of $G_a(\calC)$ by some unipotent group, and that the map $\delta : q \rightarrow p$ does not almost-split. Any such extension definably splits, but to define a splitting of the short exact sequence of binding groups, one needs a fundamental system of solutions for $q$. 

As for the necessity of (\emph{$\bigstar$}), we did not find a counterexample, but we expect there is one.

Before we end this section, we include a result about a degenerated case of splitting, which will be of use later:

\begin{lemma}\label{lem: degen-splitting}
     Let $p \in S(A)$ be a $\calC$-internal type and $f : p \rightarrow f(p)$ be an $A$-definable fibration. Suppose that $f(p)$ is weakly $\calC$-orthogonal and:
     \begin{itemize}
         \item the short exact sequence $1 \rightarrow K \rightarrow \Aut_A(p/\calC) \xrightarrow{\wt{f}} \Aut_A(f(p)/\calC) \rightarrow 1$ is $A$-definably split via some $A$-definable section $s$ and $K$ is finite,
         \item the action of $\bg{A}{f(p)}{\calC}$ satisfies \emph{($\bigstar$)}.
     \end{itemize}
    Then the groupoid $\mathrm{Gd}_A(p,f/\calC)$ is retractable, thus by Lemma \ref{lem: split-iff-retract} we get a map $\pi : p \rightarrow \pi(p)$ such that $p$ is interdefinable with $\pi(p) \times f(p)$. Moreover $\pi(p)(\calU) \subset \acl(\calC,A)$.
\end{lemma}

\begin{proof}
    By Lemma \ref{lem: ret-when-not-fun}, we know that the groupoid $\mathrm{Gd}_A(p,f/\calC)$ is retractable by the map $r$ defined by $r(f(a),f(b)) = s(\sigma)\vert_{p_{f(a)}}$, for some (any) $\sigma \in \bg{A}{f(p)}{\calC}$ with $\sigma(f(a)) = f(b)$. 

    The map $\pi$ is defined to be the quotient map of the equivalence relation $E$ given by $E(x,y)$ if and only if $r(f(a),f(b))(a) = b$, for any $a,b \models p$. By the previous paragraph, this is equivalent to $s(\sigma)\vert_{p_{f(a)}}(a) = b$. 

    Let $\sigma \in \Aut_A(\calU / \calC)$, and $a \models p$. Then there is $\tau \in K$ such that $\sigma\vert_p = s(\sigma\vert_{f(p)}) \circ \tau$. We then have:
    \begin{align*}
        \sigma(a) & = s(\sigma\vert_{f(p)}) (\tau (a)) \\
        & = s(\sigma\vert_{f(p)})\vert_{p_{f(\tau(a))}} (\tau(a)) \\
        & = r\Bigl(f\bigl(\tau(a)\bigr), f\bigl(\sigma(\tau(a))\bigr)\Bigr)(\tau(a)) \\
        & = r\Bigl(f\bigl(\tau(a)\bigr), f\bigl(\sigma(a)\bigr)\Bigr)(\tau(a)) \text{ as } \tau \in K
    \end{align*}
\noindent which implies that $E(\tau(a), \sigma(a))$. As $K$ is finite, this means that the orbit of $\pi(a)$ under $\Aut_A(\calU / \calC)$ is finite. Hence $\pi(p)(\calU) \subset \acl(\calC, A)$.
\end{proof}

Before we move on to the specific case of $f$ being the logarithmic derivative in $\mathrm{DCF}_0$, we now state and prove the strongest result we could obtain in an $\omega$-stable context, which is Theorem \ref{theo: A} of the introduction:

\begin{theorem}\label{theo: int-implies-almost-split-gen}
    Assume that $\calC$ is a pure $A$-definable algebraically closed field. Let $p \in S(A)$ be a weakly $\calC$-orthogonal type and $f : p \rightarrow q$ be an $A$-definable fibration. Assume that $q$ is $\calC$-internal with a linear nilpotent binding group and that for any $a \models p$, the fiber $p_{f(a)}$ is $\calC$-internal with diagonalizable binding group. Then $p$ is almost $\calC$-internal if and only if $f$ almost splits. 
\end{theorem}

\begin{proof}
    The right to left direction is immediate. For the left to right, we first assume that $p$ is outright $\calC$-internal. By \cite[Corollary 2.5]{jaoui2022abelian}, the binding group $\bg{A}{p}{\calC}$ is connected. There is an $A$-definable short exact sequence:
    \[1 \rightarrow K\rightarrow\bg{A}{p}{\calC} \xrightarrow{\wt{f}} \bg{A}{q}{\calC} \rightarrow 1 \]
    \noindent By \cite[Lemma 2.7]{jaoui2022abelian} $K$ is definably isomorphic to $\bg{A f(\abar)}{\tp(\abar/A f(\abar))}{\calC}$, where $\abar = a_1, \cdots , a_n$ is a Morley sequence in $p$ such that $\abar$ is a fundamental system of realizations of $p$, and $f(\abar)$ is a fundamental system of realizations of $q$. As $\bg{A f(a)}{p_{f(a)}}{\calC}$ is diagonalizable for any $f(a) \models q$, we see that $K$ is diagonalizable as a definable subgroup of the diagonalizable group $\prod\limits_{i=1}^n \bg{Af(a_i)}{\tp(a_i/Af(a_i))}{\calC}$. We would like to apply Lemma \ref{lem: nilp-splitting} to get an $A$-definable splitting of $f$. However, the group $K$ need not be connected in general, and we now show how $K$ can be replaced by its connected component by replacing $p$ with an interalgebraic type. 

    By Lemma \ref{lem: nilp-splitting}, the binding group $\bg{A}{p}{\calC}$ is definable linear nilpotent. In particular, the subgroup $K$ is central, as a subgroup of its maximal torus (which exists as $\bg{A}{p}{\calC}$ is connected). Modifying the proof of \cite[Lemma 2.2]{jaoui2022abelian}, we get that $K$ is $A$-definably isomorphic to the $\calC$-points of a diagonalizable group $\widehat{K}$, defined over $A \cap \calC$. Since $A \cap \calC$ is an algebraically closed field, the diagonalizable group $\widehat{K}$ is split over $A \cap \calC$ by \cite[8.11]{borel2012linear}, meaning it is $A \cap \calC$ definably isomorphic to an algebraic subgroup of $G_m(\calC)^k$ for some $k$. By \cite[8.7]{borel2012linear}, it is $A \cap \calC$-definably isomorphic to $\widehat{R} \times G_m(\calC)^l$, for some $l \geq 0$ and some finite, $A \cap C$-definable group $\widehat{R}$. Therefore $K$ is $A$-definably isomorphic to $R \times K^{\circ}$, for some $A$-definable finite group $R$.

    The subgroup $R$ induces an $A$-definable equivalence relation $E$ on $p$ by $E(x,y)$ if and only if there is $\sigma \in R$ such that $\sigma(x) = y$. Let $\eta$ be the quotient map, it induces an $A$-definable group morphism $\wt{\eta} : \bg{A}{p}{\calC} \rightarrow \bg{A}{\eta(p)}{\calC}$. 

    \begin{claim}
        $\ker{\wt{\eta}} = R$.
    \end{claim}

    \begin{proof}
        The kernel of $\wt{\eta}$ is the group of $\tau \in \bg{A}{p}{\calC}$ such that for all $x \models p$, there is $\sigma \in R$ such that $\tau(x) = \sigma(x)$. By Lemma \ref{lem: nilp-action-dec}, for any $x \models p$, if $\tau$ is in the unipotent radical of $\bg{A}{p}{\calC}$ and $\sigma$ is in its maximal torus, then $\tau(x) = \sigma(x)$ implies $\tau(x) = x = \sigma(x)$. As $R$ is a subgroup of the maximal torus and the action of $\bg{A}{p}{\calC}$ on $p$ is faithful, Lemma \ref{lem: nilp-action-dec} implies that $\ker(\wt{\eta})$ intersects the unipotent radical of $\bg{A}{p}{\calC}$ trivially. Again by Lemma \ref{lem: nilp-action-dec} (or rather, its proof), we see that $\ker(\wt{\eta})$ is a subgroup of the maximal torus of $\bg{A}{p}{\calC}$. Now, if $\tau \in \ker(\wt{\eta})$, there must be some $x \models p$ and $\sigma \in R$ such that $\sigma(x) = \tau(x)$. Notice that the maximal torus is central in $\bg{A}{p}{\calC}$, and $\bg{A}{p}{\calC}$ acts transitively on $p$. From this, we deduce that the maximal torus acts freely on $p$, and therefore that $\sigma = \tau$.
    \end{proof}

    Since $R < \ker(\wt{f})$, we obtain a definable map $g : \eta(p) \rightarrow f(p)$ by letting $g(\eta(x)) = f(x)$. Since $\eta$ has finite fiber, if we can show that $g : \eta(p) \rightarrow f(p)$ splits, then we automatically obtain almost splitting of $f: p \rightarrow f(p)$.

    The map $g$ induces an $A$-definable short exact sequence: 
    \[1 \rightarrow K^{\circ} \rightarrow \bg{A}{\eta(p)}{\calC} \xrightarrow{\wt{g}} \bg{A}{f(p)}{\calC} \rightarrow 1\]
    \noindent The group $\bg{A}{\eta(p)}{\calC}$ is, by Lemma \ref{lem: nilp-splitting}, a definable linear nilpotent group, i.e. definably isomorphic to the $\calC$-points of a linear nilpotent group. Moreover, the kernel $K^{\circ}$ is connected, so we can apply the second part of Lemma \ref{lem: nilp-splitting}. We only need to show that the maximal tori of $\bg{A}{\eta(p)}{\calC}$ and $\bg{A}{f(p)}{\calC}$ are $A$-definably isomorphic to the $\calC$-points of tori. Again, this can be obtained by noticing that they are central, and modifying the proof of \cite[Lemma 2.2]{jaoui2022abelian}. 

    By Lemma \ref{lem: nilp-splitting}, the short exact sequence induced by $g$ splits $A$-definably. We know that $G = \bg{A}{f(p)}{\calC} = G_u \times T$, where $G_u$ is it unipotent radical and $T$ its maximal torus. By Lemma \ref{lem: nilp-action-dec}, for any $b \models f(p)$, we have $G_u b \cap Tb = \{ b \}$. Because $G_u$ is definably unipotent, by Corollary \ref{cor: no-map-implies-star}, its action on $f(p)$ satisfies $(\bigstar)$. Again, since $\bg{A}{f(p)}{\calC}$ acts transitively on $f(p)$ and $T$ is central in it, we see that the action of $T$ on $f(p)$ is free, and therefore satisfies $(\bigstar)$. By Lemma \ref{lem: starprod} the action of $\bg{A}{f(p)}{\calC}$ on $f(p)$ satisfies $(\bigstar)$, and therefore by Lemma \ref{lem: ret-when-not-fun}, the map $g : \eta(p) \rightarrow f(p)$ splits. 

    We now explain how to reduce the case where $p$ is almost $\calC$-internal to the $\calC$-internal case, using an argument pointed out to us by the anonymous referee. Consider the map $A$-definable map $g : p \rightarrow \wt{p}$ to from $p$ its largest $\calC$-internal quotient, which is defined as follows: if $a \models p$, then $g(a)$ satisfies that
    \[\dcl(g(a)A) = \{ b \in \dcl(aA) \vert \tp(b/A) \text{ is $\calC$-internal}\} \]
    \noindent and in particular $\tp(g(a)/A)$ is $\calC$-internal. A standard argument shows that this map exists. Moreover, it must have finite fibers because $p$ is almost $\calC$-internal. As $q$ is $\calC$-internal, we obtain a map $\wt{f} : \wt{p} \rightarrow q$ such that $f = \wt{f} \circ g$. Since the fibers of $g$ are finite, almost splitting of $f$ and $\wt{f}$ are equivalent.

    The map $g$ induces, for any $a \models p$, a map on the fibers $p_{f(a)} \rightarrow \tp(g(a)/\wt{f}(a)A)$, and a routine argument shows that this induces a surjective map from the binding group of $p_{f(a)}$ to the binding group of $\tp(g(a)/\wt{f}(a)A)$. Therefore the binding group of $\tp(g(a)/\wt{f}(a)A)$ is diagonalizable as a quotient of a diagonalizable group. We apply the $\calC$-internal case of the theorem to show that the map $\wt{f}$ is almost split, which implies that the map $f$ is also almost split.
\end{proof}

\section{Logarithmic-differential pullbacks in differentially closed fields}\label{section: logdiff}

\subsection{Preliminaries on logarithmic-differential pullbacks}\label{subsection: log-diff-prel}

Using lemmas \ref{lem: fund-direct-split-retract}, \ref{lem: ret-when-not-fun} and \ref{lem: starprod}, we can use information about the splitting of definable group extensions to obtain information on the splitting of types. In this section, we are interested in the specific case, already investigated in \cite{jin2020internality} and \cite{jaoui2023relative}, of types obtained by pullback via the logarithmic derivative, in a differentially closed field of characteristic zero.

For the rest of this section, we work with differentially closed fields of characteristic zero and assume familiarity with model theory of $\mathrm{DCF}_0$, for which a reference is \cite{marker2017model}. Fix some sufficiently saturated $\calU \models \mathrm{DCF}_0$ with its derivation $\delta$. Also fix some algebraically closed differential subfield $F< \calU$. 

We let $\calC$ be the field of constants of $\calU$. It is a pure algebraically closed field, and is stably embedded. In particular, this implies that definable groups in $\calC$ are exactly algebraic groups, and definable maps between them are exactly morphisms of algebraic groups (see \cite[Chapter 7, Section 4]{marker2017model} for the equivalence of definable groups and definable morphisms with algebraic groups and morphisms).

We will consider the logarithmic derivative map:
\begin{align*}
    \dl : \calU \setminus \{ 0 \} & \rightarrow \calU \\
    x & \rightarrow \frac{\delta(x)}{x}
\end{align*}
\noindent which is a surjective morphism from $(\calU^*, \cdot )$ to $(\calU, +)$ with kernel $(\calC^*, \cdot)$.

\begin{definition}
    Let $q \in S_1(F)$. We define its pullback under the logarithmic derivative $\dl$ to be $p := \dl ^{-1}(q)$, the type of any $u$ such that $\dl(u) \models q$ and $u \not\in \acl(F, \dl(u))$.
\end{definition}

It is a known fact that this defines a unique complete type, and that each fiber $p_{\dl(u)}$ is a strongly minimal $\calC$-internal type (see \cite[Proposition 5.3]{jin2019logarithmic}). 

Consider the following three conditions on $q$, the last two of which were introduced and studied in Jin's thesis \cite{jin2019logarithmic}:

\begin{enumerate}
    \item $p = \dl^{-1}(q)$ is almost $\calC$-internal,
    \item $\dl : p \rightarrow q$ is almost split,
    \item there is an integer $k \neq 0$ such that for some $u \models p$, there are $w_1, w_2$ such that:
    \begin{itemize}
        \item $u^k = w_1 w_2$
        \item $w_1 \in \dcl(F, \dl(u))$
        \item $\dl (w_2) \in F$.
    \end{itemize}
\end{enumerate}

We will call condition $(3)$ \emph{product-splitting}.

Jin proves the implications $(3) \Rightarrow (2) \Rightarrow (1)$ in \cite[Remark 5.5]{jin2019logarithmic}. In a subsequent article \cite{jin2020internality}, Jin and Moosa consider the case where $q$ is $\calC$-internal and the generic type of an equation $\delta(x) = f(x)$, for some $f \in F(x)$ (note that it is condition $(3)$ that they call splitting). Key to their work is the following result, itself a consequence, in $\mathrm{DCF}_0$, of \cite[Chapter 3, Corollary 1.6]{marker2017model}:

\begin{fact}\label{fact: action-classification}
    Let $(G,S)$ be a definable (in $\mathrm{DCF}_0$), faithful, transitive action of a connected group $G$ on a strongly minimal set $S$. Then $(G,S)$ is definably isomorphic to either:

    \begin{enumerate}
        \item $G_a(\calC)$, $G_m(\calC)$ or $E(\calC)$ for some elliptic curve $E$ over the constants, acting regularly on itself,
        \item $G_a(\calC) \rtimes G_m(\calC)$, acting on $\calC$ by affine transformations,
        \item $\mathrm{PSL}_2(\calC)$ acting on $\mathbb{P}(\calC)$ by projective transformations.
    \end{enumerate}
\end{fact}

They prove that if the binding group of $q$ is not isomorphic to $\mathrm{PSL}_2(\calC)$, then almost internality and product splitting are equivalent, and give an almost internal, not product-split example when it is isomorphic to $\mathrm{PSL}_2(\calC)$. 

Jin and Moosa almost exclusively work with product-splitting, which has the advantage of having convenient algebraic characterizations (see Lemma \ref{lemma: jin-moosa crit}). However, splitting is more convenient to work with when using binding groups, and is key to using the methods of Section \ref{section: splitting}. It is therefore desirable to understand the connection between splitting and product splitting. In what follows, we will prove, under various additional assumptions, the equivalence of product-splitting and almost splitting, without any dimension restriction. 

The following lemma will be key in applications, and also shows that splitting of some fibrations can be reduced to splitting of the logarithmic derivative.

\begin{lemma}\label{lem: reduce to dl}
    Let $F$ be an algebraically closed field,  let $p \in S(F)$  be an almost $\calC$-internal type and $f : p \rightarrow f(p) = q$ be a fibration, with $q$ being $\calC$-internal and weakly $\calC$-orthogonal. Suppose that for some (any) $a \models p$:
    \begin{itemize}
        \item $p_{f(a)}$ is $\calC$-internal and weakly $\calC$-orthogonal,
        \item $\bg{f(a)F}{p_{f(a)}}{\calC}$ is definably isomorphic to $G_m(\calC)$.
        
    \end{itemize}
    \noindent Then there exists a type $p'$, internal to $\calC$, and a fibration $\dl : p' \rightarrow \dl(p')$ with $\dl(p') \subset \dcl(q(\calU) F)$ and the property that if $\dl : p' \rightarrow \dl(p')$ almost splits then $f : p \rightarrow f(p)$ almost splits.
\end{lemma}

\begin{proof}
    By Fact \ref{fact: almost int to int}, there is a finite-to-one definable map $\pi : p \rightarrow \pi(p)$ such that $\pi(p)$ is $\calC$-internal. In particular, for any $a \models p$, it restricts to a finite-to-one map from $p_{f(a)}$ to $\tp(\pi(a)/f(a)F)$. By Lemma \ref{lem: binding iso quotients}, there is a finite-to-one definable morphism $\bg{f(a)F}{p_{f(a)}}{\calC} \rightarrow \bg{f(a)F}{\tp(\pi(a)/f(a)F)}{\calC}$. As the former group is definably isomorphic to $G_m(\calC)$, so is the latter. 

    Fix some $a \models p$. Since $p_{f(a)}$ is weakly $\calC$-orthogonal, so is $\tp(\pi(a)/f(a)F)$. Moreover, as $q = f(p)$ is weakly orthogonal to $\calC$, we have that $\calC \cap F \langle f(a) \rangle = \calC \cap F$, which is algebraically closed. By \cite[Lemma 4.1]{jin2020internality}, there is $b \in F\langle f(a), \pi(a) \rangle \setminus F \langle f(a) \rangle ^{\mathrm{alg}}$ such that $\dl(b) \in F \langle f(a) \rangle$. We let $p' = \tp(b/F)$ and obtain a fibration $\dl : p' \rightarrow \dl(p')$ with $\calC$-internal strongly minimal fibers. Finally as $b \in \dcl(f(a), \pi(a),F)$, we see that $p'$ is $\calC$-internal. 

    We now assume that $\dl : p' \rightarrow \dl(p')$ almost splits. Then there exists $c \models r \in S(F)$, with $r$ a $\calC$-internal type, such that $b$ is interalgebraic, over $F$, with $(\dl(b),c)$ and $c \forkindep_F \dl(b)$. We are going to show that $(c, f(a))$ witnesses splitting of $f : p \rightarrow f(p)$. 

    \begin{claim}
        We have $U(c/F) = 1$.
    \end{claim}

    \begin{proof}[Proof of claim]
        We know that $b \not \in \acl(f(a),F)$, and as $\dl(b) \in \acl(f(a),F)$, we also have $b \not \in \acl(\dl(b), F)$, thus $U(\dl(b)/F) = U(b/F) - 1$. We obtain:
        \begin{align*}
            U(c/F) & = U(c/\dl(b)F) \text{ as } c \forkindep_F \dl(b) \\
            & = U(c\dl(b)/F) - U(\dl(b)/F) \\
            & = U(b/F) - U(b/F) + 1 \\
            & = 1        \end{align*}
    \end{proof}

    On the other hand, note that $\dl(b) \in F\langle f(a) \rangle$, therefore $U(b/f(a)F) = 1$, thus $U(c/f(a)F) = 1$, which implies that $c \forkindep_F f(a)$. 

    By construction $b \in \dcl(aF)$, and $c \in \acl(bF)$, so $(c, f(a)) \in \acl(aF)$. Note that $U(c,f(a)/F) = U(c/F) + U(f(a)/F) = 1 + U(f(a)/F)$. Since $\bg{f(a)F}{p_{f(a)}}{\calC}$ is definably isomorphic to $G_m(\calC)$ and $p_{f(a)}$ is weakly $\calC$-orthogonal, we see that $U(p_{f(a)}) = 1$, which implies that $U(f(a)/F) = U(a/F) - 1$, and combining this with the previous equality we get $U(c,f(a)/F) = U(a/F)$. As $(c,f(a)) \in \acl(aF)$, this implies that $a$ and $(c,f(a))$ are interalgebraic over $F$, and $f: p \rightarrow f(p)$ splits.
\end{proof}

We now start our study of pullbacks under the logarithmic derivative. The following lemma will take care of the case where the fibers of $\dl$ are not weakly $\calC$-orthogonal (and will be useful later):

\begin{lemma}\label{lem: non-w-ortho-fib implies split}
    Let $q \in S_1(F)$ be a type over any algebraically closed differential subfield of $\calU$. Let $p : = \dl ^{-1}(q)$. If for some $a \models p$, the fiber $p_{\dl(a)}$ is not weakly $\calC$-orthogonal, then there is an integer $k \neq 0$, some constant $\mu$ and some $w_1$ such that:
    \begin{itemize}
        \item $a^k = \mu w_1$
        \item $w_1 \in \dcl(F, \dl a)$.
    \end{itemize}

\noindent In particular $\dl : p \rightarrow q$ product-splits.
\end{lemma}

\begin{proof}

Let $L = F\langle \dl(a) \rangle = \dcl(F,\dl(a))$. Non-weak $\calC$-orthogonality implies that $\dcl(aL) \cap \calC \neq L \cap \calC$. As $\dl(a) \in L$, we see that $\dcl(aL) = L(a)$, the field generated by $a$ over $L$, so there are $f(X),g(X) \in L[X]$ and $e \in \calC \setminus L$ such that $\frac{f(a)}{g(a)} = e$. Write $f(x) = \sum\limits_{i=0}^n b_i x^i$ and $g(x) = \sum\limits_{j=0}^m d_j x^j$ for some $b_i, d_j \in L$ (with $b_n \neq 0$, $d_m \neq 0$ and $n \geq m$, which we may assume by replacing $e \neq 0$ by $\frac{1}{e}$). Taking a derivative, we obtain:
\[\delta(f(a))g(a)-\delta(g(a))f(a) = 0\]
\noindent which is a polynomial identity in $a$, with coefficients in $L$. Recall that $\dl^{-1}(q)$ is chosen so that $a \not\in \acl(F,\dl(a)) = \acl(L)$. Therefore all coefficients of this polynomial are zero. Computing the dominant coefficient gives:
\[ d_m\left(\delta(b
_n) + \dl(a) n b_n\right)  = b_n\left(\delta(d_m) + \dl(a)md_m\right)\]
\noindent and as $b_n \neq 0$ and $d_m \neq 0$ we get:
\[ \dl\left( a^{n-m} \right) = \dl \left( \frac{d_m}{b_n} \right) \]
\noindent and thus there is $\mu \in \calC$ such that:
\[a^{n-m} = \mu \frac{d_m}{b_n}\]
\noindent If $n \neq m$, this gives the result by picking $w_1 = \frac{d_m}{b_n}$ and $w_2 = \mu$. Otherwise, we see that $\frac{d_m}{b_m} = \mu$, for some non-zero $\mu \in \calC$. In this case, using euclidian division, we find $h \in L[X]$ with $\deg(h) < \deg(f)$ and such that $\frac{h(a)}{g(a)} = e - \mu$. As $e \not\in L$ and $\mu \in L$, we get $e-\mu \not\in L$. Since $\mu \in \calC$ and $e \in \calC$, we also get $e-\mu \in \calC$. We can then repeat the proof using these new polynomials to conclude as before. 
\end{proof}

We now prove the equivalence of almost splitting and product-splitting, as long as the fibration is not degenerate in a precise sense:

\begin{theorem}\label{theo: Jin (2) implies (3)}
    Let $q \in S_1(F)$ be a type over any algebraically closed differential subfield of $\calU$. Let $p : = \dl ^{-1}(q)$. Suppose that $p(\calU) \not\subset \acl(\dl(p)(\calU), \calC,F)$. Then the fibration $\dl:p\to q$ is almost split if and only if it is product-split, i.e. there is an integer $k \neq 0$ such that for some $a \models p$, there are $w_1,w_2$ such that:
    \begin{itemize}
        \item $a^k = w_1 w_2$
        \item $w_1 \in \dcl(F, \dl a)$
        \item $\dl (w_2) \in F$.
    \end{itemize}
\end{theorem}

\begin{proof}

As stated previously, that product-splitting implies almost splitting is essentially contained in \cite[Remark 5.5]{jin2019logarithmic}.

Suppose that the map $\dl : p \rightarrow \dl(p)$ is almost split, so there is $s \in S(F)$ that is $\calC$-internal and such that any $a \models p$ is interalgebraic with $\dl(a),b$, for some $b \models s$ with $\dl(a) \forkindep_F b$. We fix such $a$ and let $L = F \langle \dl(a) \rangle = \dcl(F, \dl(a))$, so that $b \forkindep_F L$. 

As the types $s\vert_L$ and $p_{\dl(a)}$ are interalgebraic, their binding groups must be isogenous by Corollary \ref{cor: interalg implies isogenous}. By Lemma \ref{lem: non-w-ortho-fib implies split}, we may assume that $p_{\dl(a)}$ is weakly orthogonal to $\calC$ and do so for the rest of the proof. Hence the binding group of $p_{\dl(a)}$ is definably isomorphic to $G_m(\calC)$, and the binding group of $s\vert_L$ must be as well. We want to show that $\bg{F}{s}{\calC}$ is also definably isomorphic to $G_m(\calC)$. 

First, remark that $\bg{F}{s}{\calC, \dl(p)}$ is infinite. Suppose, on the contrary, that it is finite. Then $s(\calU) \subset \acl(\calC, \dl(p)(\calU), F)$, and as $p$ is interalgebraic with $s \otimes \dl(p)$ over $F$, this would imply that $p(\calU) \subset \acl(\calC, \dl(p)(\calU),F)$, contradicting our assumption.

By Fact \ref{fact: bindgrp extension}, the binding group $\bg{L}{s\vert_L}{\calC}$ is a definable subgroup of $\bg{F}{s}{\calC}$. Moreover by Lemma \ref{lem: bind over family normal} the binding group $\bg{F}{s}{\calC, \dl(p)}$ is a normal $F$-definable subgroup of $\bg{F}{s}{\calC}$, and, by a proof similar to that of Fact \ref{fact: bindgrp extension}, is also seen to be a definable subgroup of $\bg{L}{s\vert_L}{\calC}$. Since $\bg{F}{s}{\calC, \dl(p)}$ is infinite, as $G_m(\calC)$ is connected and strongly minimal, it must be equal to $\bg{L}{s\vert_L}{\calC}$. In particular $\bg{L}{s\vert_L}{\calC}$ is normal in $\bg{F}{s}{\calC}$. 

Note that $\RM(s) = \RM(s\vert_L) = \RM(p_{\dl(a)}) = 1$, and as $F$ is algebraically closed, the type $s$ is strongly minimal. Moreover, as $s\vert_L$ is weakly $\calC$-orthogonal, forking calculus along with an automorphism argument yields that $s$ is weakly $\calC$-orthogonal, and in particular $\bg{F}{s}{\calC}$ acts transitively and faithfully on $s$. As it has a normal definable subgroup isomorphic to $G_m(\calC)$, Fact \ref{fact: action-classification} implies that $\bg{F}{s}{\calC}$ is definably isomorphic to $G_m(\calC)$.

By \cite[Lemma 4.1 (b)]{jin2020internality}, there is $e \in \dcl(bF) \setminus F $ such that $\dl (e) \in F$ (and in particular $\tp(e/F)$ is $r$-internal). Note that $e \not\forkindep_F b$, so in particular we also have $b \in \acl(eF)$, as $\RM(b/F) = 1$. So $e$ and $b$ are interalgebraic over $F$. For the rest of the proof we can therefore replace $b$ by $e$, and assume that $\dl(b) \in F$.

We have that $a$ and $b$ are non-zero, algebraically dependent over $L$, and both $\dl(a) \in L$ and $\dl(b) \in L$. By the Kolchin-Ostrowski Theorem (see \cite[top of page 1156]{kolchin1968algebraic}), we can conclude that there is some $w_1 \in L$, some $k,l \in \mathbb{Z}$ not both zero, such that $a^k b^l = w_1$. Note that $l \neq 0$, as otherwise $a \in L^{\mathrm{alg}}$, a contradiction. Therefore $k \neq 0$ and we get the conclusion by picking $w_2 = b^{-l}$.

\end{proof}

Over a field of constants, we obtain the equivalence in general :
\begin{corollary}\label{cor: prod-split-when-constant}
    Let $q \in S_1(F)$ be a type over any algebraically closed differential field of constants. Let $p : = \dl ^{-1}(q)$. Then $\dl:p\to q$ is almost split if and only if there is an integer $k \neq 0$ such that for some $a \models p$, there are $w_1,w_2$ such that:
    \begin{itemize}
        \item $a^k = w_1 w_2$
        \item $w_1 \in \dcl(F, \dl a)$
        \item $\dl (w_2) \in F$.
    \end{itemize}
\end{corollary}
\begin{proof}
    As before, we only have to prove that almost splitting implies product-splitting. In the previous proof, the only time we used the extra assumption was to show that the binding group of $s$ was definably isomorphic to $G_m(\calC)$. Keeping the notation of that proof, we now show how to obtain this if $F$ is a field of constants.
    
    By Fact \ref{fact: bindgrp extension}, the binding group $\bg{L}{s\vert_L}{\calC}$ is a definable subgroup of $\bg{F}{s}{\calC}$. Moreover, we still know that $\bg{F}{s}{\calC}$ acts transitively and faitfhfully on the strongly minimal type $s$. Over the constants, the only possible binding groups in Fact \ref{fact: action-classification} are $G_a, G_m$ and an elliptic curve (see the proof of \cite[Theorem 3.9]{freitag2022any}), and therefore the only possible case is $\bg{F}{s}{\calC}$ being definably isomorphic to $G_m(\calC)$. 
\end{proof}

We also always obtain the equivalence if the base type is orthogonal to the constants:

\begin{corollary}\label{cor: prod-split-when-ortho}
    Let $q \in S_1(F)$ be a type over any algebraically closed differential field, and assume that $q$ is orthogonal to the constants. Let $p : = \dl ^{-1}(q)$. Then $\dl:p\to q$ is almost split if and only if there is an integer $k \neq 0$ such that for some $a \models p$, there are $w_1,w_2$ such that:
    \begin{itemize}
        \item $a^k = w_1 w_2$
        \item $w_1 \in \dcl(F, \dl a)$
        \item $\dl (w_2) \in F$.
    \end{itemize}
\end{corollary}

\begin{proof}
    Again, we keep the notation of the proof of Theorem \ref{theo: Jin (2) implies (3)}. By that theorem, we may assume that $p(\calU) \subset \acl(\dl(p)(\calU), \calC , F)$, which implies that $s(\calU) \subset \acl(\dl(p)(\calU), \calC, F)$.

    As $q = \dl(p)$ is orthogonal to the constants, we have that $\bg{F}{s}{\calC, \dl(p)} = \bg{F}{s}{\calC}$. The former group is, by assumption, finite. Thus $s$ is $\calC$-algebraic. This implies, in particular, that for some (any) $a \models p$, the fiber $p_{\dl(a)}$ is not weakly $\calC$-orthogonal. By Lemma \ref{lem: non-w-ortho-fib implies split}, we get product-splitting.
\end{proof}

Finally, we prove a criteria for product-splitting of generic types of differential equations of arbitrary high order. More specifically, if $F$ is some algebraically closed differential subfield of $\calU$, we will be interested in definable sets given by solution sets of equations of the form $\delta^{m}(y) = f(y, \delta(y), \cdots , \delta^{m-1}(y))$, with $f \in F(x_0, \cdots , x_{m-1})$. Note that we can write $f = \frac{P}{Q}$ with $P,Q \in F[x_0, \cdots , x_{m-1}]$ without common factors, and that our definable set is then given by solutions of the equation $\delta^{m}(y)Q(y, \delta(y), \cdots , \delta^{m-1}(y)) - P(y, \delta(y), \cdots , \delta^{m-1}(y)) = 0$. As the polynomial $x_m P(x_0, \cdots , x_{m-1}) - Q(x_0, \cdots, x_{m-1}) \in F[x_0, \cdots , x_m]$ is irreducible, this definable set has a unique generic type.

Recall that if $(F,\delta)$ is any differential field, then the field $F(x_0, \cdots , x_m)$ can be equipped with the derivations $\frac{\partial }{\partial x_i}$ with respect to $x_i$, as well as the derivation $\delta^F$, which treats all $x_i$ as constants and equals $\delta$ on $F$. If $P \in F[x_0, \cdots , x_m]$, then $\delta^F(P) = P^{\delta}$, i.e. the polynomial obtained by differentiating the coefficients of $P$. Using this notation, we can compute, for any tuple $(a_0, \cdots , a_{m-1}) \in \calU^m$ and $h \in F(x_0, \cdots , x_{m-1})$, that $\delta(h(\overline{a})) = \sum\limits_{i = 0}^{m-1} \frac{\partial h}{\partial x_i} \delta(a_i) + \delta^F(h)(\overline{a})$. We obtain the following:

\begin{lemma}\label{lemma: jin-moosa crit}
    Let $F$ be an algebraically closed differential field, a rational function $f \in F(x_0, \cdots , x_{m-1})$ and $q$ the generic type of $\delta^{m}(y) = f(y, \delta(y), \cdots , \delta^{m-1}(y))$. Let $p := \dl^{-1}(q)$, then $\dl : p \rightarrow q$ is product-split if and only if there are a non-zero $h \in F(x_0, \cdots , x_{m-1})$, some $e \in F$ and some integer $k \neq 0$ such that 
    \[ (kx_0-e)h = \sum\limits_{i=0}^{m-2} \frac{\partial h}{\partial x_i}x_{i+1} + \frac{\partial h}{\partial x_{m-1}}f + \delta^F(h) \text{ .} \] 
    
\end{lemma}
\begin{proof}
    Suppose that there are such $h,e$ and $k$, consider some $u \models p$ and $a = \dl(u) \models q$. Also denote $(a, \delta(a), \cdots , \delta^{m-1}(a)) = \overline{a}$. As $a$ realizes the generic type of $\delta^{m}(y) = f(y, \delta(y), \cdots , \delta^{m-1}(y))$, which is of order $m$, we have that $h(\overline{a})$ is well-defined and non-zero, as otherwise $a$ would satisfy a differential equation of order $m-1$.
    
    Let $w_1 = h(\overline{a})$ and $w_2 = \frac{u^k}{h(\overline{a})}$, so that $u^k = w_1 w_2$, we need to show that $w_1 \in \dcl(F, \dl(u))$ and $\dl(w_2) \in F$. The first part is immediate, as $a = \dl(u) $ and $h \in F(x_0, \cdots , x_{m-1})$. We also compute:

    \begin{align*}
        \dl(w_2) & = ka - \dl(h(\overline{a})) \\
        & = ka - \frac{\delta(h(\overline{a}))}{h(\overline{a})} \\
        & = ka - \frac{\sum\limits_{i=0}^{m-1} \frac{\partial h}{\partial x_i}(\overline{a}) \delta^{i+1}(a) + \delta^F(h)(\overline{a})}{h(\overline{a})} \\
        & = ka - \frac{\sum\limits_{i=0}^{m-2} \frac{\partial h}{\partial x_i}(\overline{a}) \delta^{i+1}(a) + \frac{\partial h}{\partial x_{m-1}}(\overline{a}) f(\overline{a}) + \delta^F(h)(\overline{a})}{h(\overline{a})}  \\
        & = ka - (ka-e) \\
        & = e \in F \text{ .}
    \end{align*}

    For the converse, suppose that $\dl : p \rightarrow q$ is product-split. There are $u \models p$, some integer $k \neq 0$ and $w_1, w_2$ such that:
    \begin{itemize}
        \item $u^k = w_1 w_2$
        \item $w_1 \in \dcl(F, \dl u)$
        \item $\dl (w_2) \in F$.
    \end{itemize}
    Note that neither $w_1$ or $w_2$ is zero as $u \neq 0$. We let $a = \dl(u)$ and $e = \dl(w_2)$. Again denote $\overline{a} = (a, \delta(a), \cdots , \delta^{m-1}(a))$. As $\delta^m(a) = f(\overline{a})$ and $w_1 \in \dcl(F,a)$, there is a non-zero $h \in F(x_0, \cdots , x_{m-1})$ such that $w_1 = h(\overline{a})$. Applying the logarithmic derivative to the equality $u^k = w_1 w_2$, we compute:

    \begin{align*}
        (ka-e)h(\overline{a}) & = \delta(h(\overline{a})) \\
        & = \sum\limits_{i=0}^{m-1} \frac{\partial h}{\partial x_i}(\overline{a}) \delta^{i+1}(a) + \delta^F(h)(\overline{a}) \\
        & = \sum\limits_{i=0}^{m-2} \frac{\partial h}{\partial x_i}(\overline{a}) \delta^{i+1}(a) + \frac{\partial h}{\partial x_{m-1}}(\overline{a}) f(\overline{a}) + \delta^F(h)(\overline{a}) \text{ .}
    \end{align*}
    \noindent Because $a$ realizes the generic type of $\delta^{m}(y) = f(y, \delta(y), \cdots , \delta^{m-1}(y))$ and this equation is of order $m-1$, this implies:

    \[(kx_0-e)h = \sum\limits_{i=0}^{m-2} \frac{\partial h}{\partial x_i}x_{i+1} + \frac{\partial h}{\partial x_{m-1}}f + \delta^F(h) \text{ .}\]
\end{proof}

Under the additional assumption of Theorem \ref{theo: Jin (2) implies (3)}, we get a criteria for almost splitting:

\begin{corollary}\label{cor: almost-split-char}
    Let $F$ be an algebraically closed differential field, consider some $f \in F(x_0, \cdots , x_{m-1})$ and $q$ the generic type of $\delta^{m}(y) = f(y, \delta(y), \cdots , \delta^{m-1}(y))$. Let $p := \dl^{-1}(q)$ and assume that $p(\calU) \not\subset \acl(\dl(p)(\calU),\calC,F)$. Then $\dl : p \rightarrow q$ is almost split if and only if there are a non-zero $h \in F(x_0, \cdots , x_{m-1})$, some $e \in F$ and some integer $k \neq 0$ such that 
    \[ (kx_0-e)h = \sum\limits_{i=0}^{m-2} \frac{\partial h}{\partial x_i}x_{i+1} + \frac{\partial h}{\partial x_{m-1}}f + \delta^F(h) \text{ .} \] 
\end{corollary}

\begin{proof}
    Immediate consequence of Theorem \ref{theo: Jin (2) implies (3)} and Lemma \ref{cor: jin-moosa crit}.
\end{proof}

\subsection{Splitting when the binding group is nilpotent}\label{subsection: nilpotent}

In this subsection, we prove that for any $\calC$-internal and weakly $\calC$-orthogonal type $q$ over any algebraically closed field $F$, if $q$ has nilpotent binding group, then $\dl^{-1}(q)$ is almost $\calC$-internal if and only if $\dl: \dl^{-1}(q) \rightarrow q$ splits. This can be seen as a partial generalization of Theorem A of \cite{jin2020internality}. 

We will sometimes return to the notation of Subsection \ref{subsection: jordan}: for example, by a definable nilpotent linear group, we mean a definable group definably isomorphic to the $\calC$-points of a nilpotent linear algebraic group.

\begin{theorem}\label{theo: splits if nilpotent base}
    Let $F$ be an algebraically closed differential field and let $q \in S_1(F)$ be a $\calC$-internal, weakly $\calC$-orthogonal type with binding group definably isomorphic to the $\calC$-points of a nilpotent linear algebraic group. Then the following are equivalent:

    \begin{enumerate}
        \item $\dl^{-1}(q)$ is almost $\calC$-internal,
        \item $\dl: \dl^{-1}(q) \rightarrow q$ product-splits,
        \item $\dl: \dl^{-1}(q) \rightarrow q$ almost splits.
    \end{enumerate}

\end{theorem}

\begin{proof}

    The implication $(2) \Rightarrow (3)$, as previously discussed, is essentially given by \cite[Remark 5.5]{jin2019logarithmic}, and the implication $(3) \Rightarrow (1)$ is immediate. Therefore we only prove $(1) \Rightarrow (2)$. In fact, we will prove $(1) \Rightarrow (3)$ and $\bigl( (1) \text{ and }(3) \bigr) \Rightarrow (2)$. 
    
    Assume that $p = \dl^{-1}(q)$ is almost $\calC$-internal. We will start by showing that this implies almost splitting of $\dl$. To use binding groups, we need to reduce to the case where $p$ is $\calC$-internal. 

    First note that if $p_{\dl(a)}$ is not weakly $\calC$-orthogonal, by Lemma \ref{lem: non-w-ortho-fib implies split} $\dl: \dl^{-1}(q) \rightarrow q$ product-splits. So we now assume that $p_{\dl(a)}$ is weakly $\calC$-orthogonal. As $q$ is also weakly $\calC$-orthogonal, we see that $p$ is weakly $\calC$-orthogonal.
    
    By Lemma \ref{lem: reduce to dl}, there is another type $p'$, internal to $\calC$, and a fibration $\dl : p' \rightarrow \dl(p')$ with $\dl(p') \subset \dcl(q(\calU)F)$, with the property that if $\dl : p' \rightarrow \dl(p')$ almost splits then $\dl : p \rightarrow \dl(p)$ almost splits. Because $p$ is weakly $\calC$-orthogonal, the proof of Lemma \ref{lem: reduce to dl} shows that $p'$ also is. 
    
    We see that $\dl(p')$ is weakly $\calC$-orthogonal, and its binding group also is definably isomorphic to the $\calC$-points of a nilpotent linear algebraic group. Moreover if $\dl: p' \rightarrow \dl(p')$ almost splits, so does $\dl : p \rightarrow \dl(p)$. Therefore, replacing $p$ by $p'$, we may assume that $p$ is $\calC$-internal. By Theorem \ref{theo: int-implies-almost-split-gen}, the fibration $\dl: p \rightarrow q$ is almost split.    

    We still need to obtain product-splitting. By Theorem \ref{theo: Jin (2) implies (3)}, we may assume that $p(\calU) \subset \acl(\dl(p)(\calU), \calC , F)$. Using similar ideas as in the proof of Corollary \ref{cor: interalg implies isogenous}, we see that this implies that $\bg{F}{p}{\calC}$ is finite. This group is the kernel of $\dl : p \rightarrow \dl(p)$. By Lemma \ref{lem: degen-splitting}, we obtain that the fibers of $f$ are not weakly $\calC$-orthogonal, which is a contradiction, as we assumed that they were.

\end{proof}

We deduce some internality criteria:

\begin{corollary}\label{cor: jin-moosa crit}
    Let $F$ be an algebraically closed differential field, some rational function $f \in F(x_0, \cdots , x_{m-1})$, and consider $q$, the generic type of $\delta^{m}(y) = f(y, \delta(y), \cdots , \delta^{m-1}(y))$. If $q$ is $\calC$-internal, weakly $\calC$-orthogonal and has a nilpotent binding group, then $p := \dl^{-1}(q)$ is almost $\calC$-internal if and only if there are a non-zero $h \in F(x_0, \cdots , x_{m-1})$, some $e \in F$ and some integer $k \neq 0$ such that 
    \[ (kx_0-e)h = \sum\limits_{i=0}^{m-2} \frac{\partial h}{\partial x_i}x_{i+1} + \frac{\partial h}{\partial x_{m-1}}f + \delta^F(h)\text{ .}\]   
\end{corollary}

\begin{proof}

    This is an immediate consequence of Lemma \ref{lemma: jin-moosa crit} and Theorem \ref{theo: splits if nilpotent base}.  
\end{proof}

In practice, we can now sometimes reduce the question of internality of a logarithmic differential pullback to some valuation computation. Consider some differential field $F$ and the rational function field $F(x_0, \cdots , x_{m-1})$, recall that it is equipped with the derivations $\frac{\partial}{\partial x_i}$, as well as the derivation $\delta^F$ that is equal to $\delta$ on $F$, and $\delta^F(x_i) = 0$ for all $i$.

We can also equip it with valuations: see $F(x_0, \cdots  , x_{m})$ as $F(x_1, \cdots , x_{m})(x_0)$, and equip it with the valuation given, for any $P,Q \in F(x_1, \cdots , x_{m})[x_0]$, by $v(P) = - \mathrm{deg}(P)$, and $v(\frac{P}{Q}) = v(P) - v(Q)$. We have, for all $P,Q \in F(x_1, \cdots , x_{m})[x_0]$ and $i$, that $v\left(\frac{\partial P}{\partial x_i}\right) \geq v(P)$ and $v(\delta^F(P)) \geq v(P)$, from which we deduce that $v\left(\frac{\partial \frac{P}{Q}}{\partial x_i}\right) \geq v\left( \frac{P}{Q} \right) $ and $v(\delta^F(\frac{P}{Q})) \geq v(\frac{P}{Q})$. We could do this for any $0 \leq i \leq m$, and we denote $v_i$ the valuation obtained.

As an example of application, we look at linear differential equations:

\begin{corollary}
    Let $F$ be an algebraically closed differential field, and consider $L(x_0, \cdots, x_{m-1}) = \sum\limits_{i=0}^{m-1} c_i x_i + c \in F[x_0, \cdots , x_{m-1}]$. Let $q$ be the generic type of $\delta^m(x) = L(\delta^{m-1}(x) , \cdots , \delta(x), x)$, which is always $\calC$-internal. If $q$ is weakly $\calC$-orthogonal and has nilpotent binding group, then $\dl^{-1}(q)$ is not almost $\calC$-internal.
\end{corollary}

\begin{proof}
    Suppose, for a contradiction, that $\dl^{-1}(q)$ is almost $\calC$-internal. By Corollary \ref{cor: jin-moosa crit}, there are a non-zero $h \in F(x_0, \cdots, x_m)$, some $e \in F$ and some integer $k \neq 0$ such that we have the following equality, in $F(x_0, \cdots , x_{m-1})$:
    \[ (kx_0-e)h = \sum\limits_{i=0}^{m-2} \frac{\partial h}{\partial x_i}x_{i+1} + \frac{\partial h}{\partial x_{m-1}}\left( \sum\limits_{j=0}^{m-1}c_j x_ j + c \right) + \delta^F(h) \text{ .}\] 
    Let $v_0$ be the valuation in $F(x_0, \cdots , x_{m-1})$ with respect to $x_0$. We can write $h = \wt{h}x_0^{-v_0(h)} + g(x_0, \cdots , x_{m-1})$ with $\wt{h} \in F(x_1, \cdots , x_{m-1})$ and $v_0(g) > v_0(h)$. We have:
    \[ k\wt{h} = c_0 \frac{\partial \wt{h}}{\partial x_{m-1}} \text{ .}\]
    \noindent Now consider the valuation $v_{m-1}$ with respect to $x_{m-1}$, the left-hand side has valuation $v_{m-1}(\wt{h})$, but the right-hand side has valuation strictly greater than $v_{m-1}(\wt{h})$, unless $\wt{h} = 0$, which is a contradiction, as then $g = h$, but $v_0(g)> v_0(h)$ (note that this also works in the case where $m = 1$ and $\wt{h} \in F$).
\end{proof}

If $F$ is a field of constants, then this would apply as long as the generic type of $\delta^m(x) = L(\delta^{m-1}(x), \cdots , \delta(x),x)$ is weakly $\calC$-orthogonal, as its binding group is always abelian in that case (see the proof of \cite[Theorem 3.9]{freitag2022any} for a justification). Note that this is not always the case: for example, the generic type of $\delta^{m}(x) = 0$ is never weakly $\calC$-orthogonal as long as $m > 1$.

Over non-constant parameters, we could pick some differential transcendental $c$, and let $F = \mathbb{Q}\langle c \rangle ^{\mathrm{alg}}$. Then this would apply to the generic type of $\delta^m(x) = c$, which one can easily show is $\calC$-internal, weakly $\calC$-orthogonal and with unipotent (but non-abelian) binding group.

\subsection{Non-splitting when the binding group is an elliptic curve}\label{subsection: elliptic}

We now want to use our methods to exhibit more uniformly internal maps that do not split. In \cite[bottom of page 5]{jin2020internality}, Jin and Moosa ask whether there exists a strongly minimal type $q \in S_1(F)$, with $F$ algebraically closed, such that $q$ is $\calC$-internal with binding group an elliptic curve, with the type $\dl^{-1}(q)$ almost $\calC$-internal, but the map $\dl : \dl^{-1}(q) \rightarrow q$ not almost split. Using Lemma \ref{lem: fund-direct-split-retract}, we will show that there exist many such types. 

To deal with algebraic closure issues, we will need to use the machinery  of the $\ext(\cdot, \cdot)$ functor classifying extensions of commutative algebraic groups. A good reference is \cite[Chapter 7]{serre2012algebraic}, which we will use in the rest of this section. We now recall some of the exposition and facts found in that chapter.

Fix some base algebraically closed field $k$, over which everything will be defined, and a $\vert k \vert$-saturated algebraically closed field $\calC \supset k$ (it will be the constants in our application) in which all definable sets and algebraic groups will live. Recall that an extension of algebraic groups is a short exact sequence:
\[1 \rightarrow B \rightarrow C \rightarrow A \rightarrow 1\]
\noindent of algebraic groups, where the morphisms are maps of algebraic groups. Two such extensions $C$ and $C'$ are isomorphic if there exists a map $f : C \rightarrow C'$ making the diagram:
\begin{center}
\begin{tikzcd}
    1 \arrow[r] & B \arrow[r] \arrow[d, "\id"] & C \arrow[r] \arrow[d, "f"] & A \arrow[r] \arrow[d, "\id"] & 1 \\
    1 \arrow[r] & B \arrow[r] & C' \arrow[r] & A \arrow[r] & 1 
\end{tikzcd}
\end{center}
\noindent commute. In that case, $f$ must be an isomorphism. For any two commutative algebraic groups $A$ and $B$, denote $\ext(A,B)$ the set of isomorphism classes of commutative group extensions. We will now work with commutative groups and switch to additive notation. Here is a summary of the material of \cite[VII.1]{serre2012algebraic} making $\ext(\cdot , \cdot)$ into a functor.

Given a morphism $f : B \rightarrow B'$, there exists a unique extension $C' \in \ext(A,B')$ and a map $F : C \rightarrow C'$ making the following diagram:
\begin{center}
\begin{tikzcd}
    0 \arrow[r] & B \arrow[r] \arrow[d, "f"] & C \arrow[r] \arrow[d, "F"] & A \arrow[r] \arrow[d, "\id"] & 0 \\
    0 \arrow[r] & B' \arrow[r] & C' \arrow[r] & A \arrow[r] & 0 
\end{tikzcd}
\end{center}
\noindent commute. We denote this extension $f_*(C)$. If $g :  A' \rightarrow A$ is a morphism, there is a similar construction of $g^*(C)$. We give details on the construction, as we will need them later. Let:
\[0 \rightarrow B \rightarrow C \xrightarrow{f} A \rightarrow 0\]
\noindent be an element of $\ext (A,B)$. There is a unique $C'  \in \ext (A',B)$ and map $G : C' \rightarrow C$ making the diagram:
\begin{center}
\begin{tikzcd}
    0 \arrow[r] & B \arrow[r] \arrow[d, "\id"] & C' \arrow[r] \arrow[d, "G"] & A' \arrow[r] \arrow[d,"g"] & 0 \\
    0 \arrow[r] & B \arrow[r] & C \arrow[r, "f"] & A \arrow[r] & 0 
\end{tikzcd}
\end{center}
\noindent commute. Moreover, the group $C'$ is the subgroup of $A' \times C$ consisting of pairs $(a,c)$ such that $g(a) = f(c)$, and the maps $C' \rightarrow A'$ and $G : C' \rightarrow C$ are the natural projections.

The following summarizes the properties we will need (see \cite[Chapter 7, 1.1]{serre2012algebraic}):

\begin{fact}\label{fact: ext-is-functor}

Let $\mathcal{A}$ be the category of commutative algebraic groups. The previous construction makes $\ext(\cdot, \cdot)$ into an additive bifunctor on $\mathcal{A} \times \mathcal{A}$, contravariant in the first coordinate and covariant in the second. In particular $\ext(A,B)$ is always an algebraic group, and the maps $f_*$ and $g^*$ are morphisms of algebraic groups. 
    
\end{fact}

We will also need to identify the specific algebraic group $\ext(E,G_m)$, where $E$ is an elliptic curve. Recall that for any abelian variety $A$, we can form its dual abelian variety $\wt{A}$, which parametrizes the topologically trivial line bundles on $A$ (see \cite[II.8]{mumford1974abelian}). The dual of an elliptic curve is also an elliptic curve. We then have the following (see \cite[Chapter 7, 3.16]{serre2012algebraic}):

\begin{fact}\label{fact: ext(E)}
    The group $\ext(E,G_m)$ is isomorphic to the dual elliptic curve $\widehat{E}$.
\end{fact}

Finally, if $f : G_m \rightarrow G_m$ is a non-trivial morphism of algebraic groups, we will need information on the kernel of the map $f_* : \ext (E,G_m) \rightarrow \ext (E,G_m)$. We include the proof for the reader's (and authors') comfort, even though it may follow immediately from known homological algebra facts.

\begin{proposition}\label{prop: ext-map-fin-ker}
    Let $E$ be an elliptic curve, and $g : G_m \rightarrow G_m$ be a non-trivial morphism of algebraic groups. Then $g_* : \ext (E,G_m) \rightarrow \ext (E,G_m)$ has finite kernel.
\end{proposition}

\begin{proof}
    The map $g : G_m \rightarrow G_m$ has finite kernel, denote it $H$. By \cite[Chapter 7, 4.23, Theorem 12]{serre2012algebraic}, the functor $\ext (E, \cdot)$ is exact on the category of linear algebraic groups, and we thus obtain a short exact sequence:
    \[0 \rightarrow \ext(E,H) \rightarrow \ext (E,G_m) \xrightarrow{g_*} \ext (E,G_m) \rightarrow 0\]
    \noindent It is therefore enough to show that $\ext (E, H)$ is finite, whenever $H$ is a finite algebraic group. 

    To do so, we show that there is a surjective morphism from $\hom (\prescript{}{n}{E},H)$ to $\ext (E,H)$, where $\prescript{}{n}{E}$ is the $n$-torsion of $E$. This torsion group is well-known to be finite (see \cite[II.4]{mumford1974abelian}), so $\hom (\prescript{}{n}{E},H)$ is finite as well. 

    Let $n \in \mathbb{N}$, we have maps $\nu_n : H \rightarrow H$ and $\mu_n : E \rightarrow E$ given by sending $x$ to $nx$. It is well known (see \cite[II.4]{mumford1974abelian}) that $\mu_n$ is surjective and has finite kernel. Pick $n$ to be a multiple of the order of $H$, so that $\nu_n$ is the zero map, or in other words, so that $nh =  0$ for all $h \in H$. We have a short exact sequence:
    \[ 0 \rightarrow \prescript{}{n}{E} \xrightarrow{\iota} E \xrightarrow{\mu_n } E \rightarrow 0 \]
    \noindent which is an element $e$ of $\ext (E, \prescript{}{n}{E})$. Any $\phi \in \hom (\prescript{}{n}{E}, G)$ gives rise to a map $\phi_* : \ext (E, \prescript{}{n}{E}) \rightarrow \ext (E,G)$, and in particular we obtain $\phi_* (e) \in \ext (E,G)$. To summarize, we have constructed a map:
    \begin{align*}
        d : \hom (\prescript{}{n}{E} , G) & \rightarrow \ext (E , G) \\
        \phi & \rightarrow \phi_*(e)
    \end{align*}
    \noindent By \cite[Chapter 7, 1.2, Proposition 2]{serre2012algebraic}, we obtain an exact sequence: 
    \begin{center}
    \begin{tikzcd}
        0 \arrow[r] & \hom (E,H) \arrow[r, "\cdot \circ \mu_n"] & \hom (E,H) \arrow[r, "\cdot \circ \iota"] & \hom(\prescript{}{n}{E}, H) \arrow[lld, "d"'] \\
         & \ext (E,H) \arrow[r, "\mu_n^*"'] & \ext(E,H) \arrow[r, "\iota^*"'] & \ext(\prescript{}{n}{E}, H)
    \end{tikzcd}
    \end{center}
    
     \noindent We show that $d$ is surjective by proving that the map $\mu_n^*$ is the zero map. Let $0 \rightarrow H \rightarrow C \xrightarrow{f} E \rightarrow 0$ be an element of $\ext (E,H)$. We consider the group $C' < E \times C$ consisting of pairs $(e,c)$ such that $ne = f(c)$ and obtain the extension $\mu_n^* (C)$ as the top row of the commutative diagram:

    \begin{center}
    \begin{tikzcd}
    0 \arrow[r] & H \arrow[r] \arrow[d, "\id"] &  C' \arrow[r,"f'"] \arrow[d, "G"] & E \arrow[r] \arrow[d,"\mu_n"] & 0 \\
    0 \arrow[r] & H \arrow[r] & C \arrow[r, "f"] & E \arrow[r] & 0 
\end{tikzcd}
\end{center}

\noindent We need to show that the top short exact sequence is the trivial extension, or equivalently, that it splits (definably and without needing extra parameters). 

We show that for any $e \in E$, there is a unique $x \in nC'$ such that $f'(x) = e$, this will define a section. We know that $\mu_n$ and $f$ are surjective, so there are $\wt{e} \in E$ and $c \in C$ such that $n\wt{e} = e = f(c)$. Then $(\wt{e},c) \in C'$ and $f'(n(\wt{e},c)) = e$, giving existence. For uniqueness, suppose that $f'(n(e_1,c_1)) = f'(n(e_2,c_2))$. This implies that $ne_1 = ne_2$, and as $(e_1,c_1), (e_2,c_2) \in C'$, that $f(c_1) = f(c_2)$. Hence $c_1 - c_2 \in H$, and by choice of $n$, we obtain that $nc_1 = nc_2$. Therefore $n(e_1, c_1) = n(e_2,c_2)$. 

Thus we have obtained a section to $f'$, which is immediately seen to be definable without extra parameters. That it is a morphism is left to the reader.
    
\end{proof}

To complete our proof, we will need to exhibit a $\calC$-internal type with binding group isomorphic to a semiabelian variety. This follows from Kolchin's solution to the inverse Galois problem for strongly normal extensions (see \cite[Theorem 2]{kolchin1955galois} and also \cite[Proposition 15.1]{kovacic2006geometric}), as well as some translation to model-theoretic language. We give some details. Our exposition closely follows Marker's in \cite[Chapter 2, Section 9]{marker2017model}. 

Let $K,L$ be differential subfields of $\calU$, and denote $\calC_K = \calC \cap K$ and $\calC_L = \calC \cap L$. We say that $K < L$ is a \emph{strongly normal} extension if:

\begin{enumerate}
    \item $\calC_K = \calC_L$ is algebraically closed,
    \item $L/K$ is finitely generated,
    \item for any $\sigma \in \Aut_K(\calU)$, we have that $\langle L , \calC \rangle = \langle \sigma(L), \calC \rangle$.
\end{enumerate}

The \emph{differential Galois group} $G(L/K)$ is the group of differential automorphisms of $L$ fixing $K$. The \emph{full} differential Galois group $\mathrm{Gal}(L/K)$ is the group of differential automorphisms of $\langle L, \calC \rangle$ fixing $\langle K, \calC \rangle$. 

Since $L/K$ is strongly normal, there is some tuple $\overline{a}$ such that $L = K \langle \overline{a} \rangle$. Let $p = \tp(\overline{a}/K)$. Marker's proof of \cite[Theorem 9.5]{marker2017model} (and the discussion following it) shows that $\tp(\overline{a}/K)$ is $\calC$-internal, and that there is an algebraic group $G$ such that $G(L/K)$ (resp. $\mathrm{Gal}(L/K)$) is isomorphic to $G(\calC_K)$ (resp. $G(\calC)$). 

It is easy to see that $\mathrm{Gal}(L/K)$ is definably isomorphic to the binding group $\bg{K}{p}{\calC}$. In other words, the binding group of $p$ is isomorphic to the $\calC$-points of the Galois group of the strongly normal extension $L/K$. By Kolchin's solution to the inverse Galois problem, any connected algebraic group is the Galois group of some strongly normal extension, over some differential field $F$, which we can assume to be algebraically closed. We can then take the type $p$ obtained previously, which is weakly $\calC$-orthogonal as $\calC_K = \calC_L$. We have obtained:

\begin{fact}\label{fact: inverse galois}
    Fix a connected algebraic group $G$ defined over $\calC$. There exists an algebraically closed differential field $F$ and some $p \in S(F)$ such that $p$ is $\calC$-internal, weakly $\calC$-orthogonal, and $\bg{F}{p}{\calC}$ is definably isomorphic to $G(\calC)$.
\end{fact}

Finally, we will use the Galois correspondence for binding groups. We refer the reader to \cite[Theorem 2.3]{sanchez2017some} for a proof, as well as the closest account to what we need that we could find. From that theorem, we deduce the following:

\begin{fact}
    Let $F$ be an algebraically closed field, and $p \in S(F)$ a $\calC$-internal, fundamental and weakly $\calC$-orthogonal type. Fix some $a \models p$. If $H$ is an $F$-definable normal subgroup of $\bg{F}{p}{\calC}$, then there is $b \in \dcl(aF)$ such that:
    \begin{itemize}
        \item $H = \{ \sigma \in \bg{F}{p}{\calC}, \text{ some (any) lift of } \sigma \text{ to } \calU \text{ fixes } b\}$,
        \item $\tp(b/F)$ is $\calC$-internal, weakly $\calC$-orthogonal and fundamental, and the binding group $\bg{F}{\tp(b/F)}{\calC}$ is definably isomorphic to $\bg{F}{p}{\calC}/H$,
        \item $\tp(a/bF)$ is $\calC$-internal, weakly $\calC$-orthogonal and fundamental, and the binding group $\bg{F}{\tp(a/bF)}{\calC}$ is definably isomorphic to $H$.
    \end{itemize}
    If $H$ is connected, we also obtain that $\tp(a/bF)$ is stationary (see \cite[Chapter 1, Lemma 6.16]{pillay1996geometric}).        

\end{fact}

From this we obtain:

\begin{fact}\label{fact: galois correspondence}
    Let $F$ be an algebraically closed field, let $p \in S(F)$ be a $\calC$-internal, fundamental and weakly $\calC$-orthogonal type. If $H$ is a connected $F$-definable normal subgroup of $\bg{F}{p}{\calC}$, then there is an $F$-definable fibration $f : p \rightarrow f(p)$ with $\calC$-internal fibers, such that $f(p)$ and $p_{f(a)}$ (for any $f(a) \models f(p)$) are $\calC$-internal and fundamental, with $\bg{F}{p_{f(a)}}{\calC}$ definably isomorphic to $H$, and giving rise to the following short exact sequence:
    \[1 \rightarrow H \rightarrow \bg{F}{p}{\calC} \xrightarrow{\wt{f}} \bg{F}{f(p)}{\calC} \rightarrow 1 \text{ .}\]
    \noindent
\end{fact}

We are now equipped to prove the following:

\begin{theorem}\label{theo: non-split-pullback}
    There exist an algebraically closed field $F$, some type $p \in S(F)$ that is $\calC$-internal, fundamental and weakly $\calC$-orthogonal, and a definable fibration $f: p \rightarrow f(p)$, such that:
    \begin{itemize}
        \item $\bg{F}{p}{\calC}$ is definably isomorphic to $D(\calC)$, where $D$ is a non-definably split extension of an elliptic curve by $G_m$,
        \item $f$ does not almost split.
        \item for any $a \models p$, the type $p_{f(a)}$ is $\calC$-internal, weakly $\calC$ orthogonal, fundamental, with binding group definably isomorphic to $G_m(\calC)$.
    \end{itemize}
\end{theorem}

\begin{proof}
    By Fact \ref{fact: ext(E)}, there is an algebraic group $D$ that is a non-definably split extension of an elliptic curve $E$ by $G_m$, i.e. we have a map $g : D \rightarrow E$ with kernel $G_m$. Moreover, we can assume that $D$ does not belong to any finite subgroup of the elliptic curve $\mathrm{Ext}(E,G_m)$ (equivalently $D$ is not a torsion point).
    
    By Fact \ref{fact: inverse galois}, there is an algebraically closed field $F$ and $p \in S(F)$ such that $p$ is $\calC$-internal and weakly $\calC$-orthogonal, with $\bg{F}{p}{\calC}$ definably isomorphic to $D(\calC)$, where $G$ is a non-definably split extension of an elliptic curve by $G_m$. As was observed in \cite[Fact 2.4]{jaoui2022abelian}, we may assume that $p$ is fundamental by taking a Morley power. 

    In particular the group $\bg{F}{p}{\calC}$ has a definable Chevalley decomposition in the sense of \cite[Fact 2.8]{jaoui2022abelian}, meaning a subgroup maximal among definable subgroups $F$-definably isomorphic to the $\calC$-points of a linear algebraic group. This subgroup $L$ is called the \emph{linear part} of $\bg{F}{p}{\calC}$. It is easy to show, in this case, that the linear part of $\bg{F}{p}{\calC}$ is the kernel of the map induced by $g$ on $\bg{F}{p}{\calC}$ (in particular it is definably isomorphic to $G_m(\calC)$). Moreover, the linear part is $F$-definable and normal by \cite[Fact 2.8]{jaoui2022abelian}. 

    By Fact \ref{fact: galois correspondence}, there is an $F$-definable fibration $f : p \rightarrow f(p)$ giving rise to the short exact sequence:
    \[ 1 \rightarrow L \rightarrow \bg{F}{p}{\calC} \xrightarrow{\wt{f}} \bg{F}{f(p)}{\calC} \rightarrow 1\]
    \noindent and $\bg{F}{f(p)}{\calC}$ is definably isomorphic to $E(\calC)$. 

    Suppose, for a contradiction, that $f$ does almost split. Then there exist, by Proposition \ref{pro: reduce to int}, some $\calC$ internal type $s \in S(F)$, some $a \models p$  and $b \models s$ such that $a$ is interalgebraic with $(f(a),b)$ over $F$, with $b \forkindep_F f(a)$. 
    
    Since both $p$ and $f(p)$ are fundamental, by \cite[Lemma 2.7]{jaoui2022abelian}  $\bg{f(a)F}{p_{f(a)}}{\calC}$ is definably isomorphic to $L$, and hence to $G_m(\calC)$. By Lemma \ref{lem: binding iso quotients} the binding groups $\bg{f(a)F}{p_{f(a)}}{\calC}$ and $\bg{f(a)F}{\tp(b/f(a)F)}{\calC}$ are isogenous, thus $\bg{f(a)\calC}{\tp(b/f(a)F)}{\calC}$ is isogenous to $G_m(\calC)$. Also $\bg{f(a)F}{\tp(b/f(a)F)}{\calC}$ is a definable subgroup of $\bg{F}{s}{\calC}$ by Fact \ref{fact: bindgrp extension}. 

    As $L$ acts regularly on $p_{f(a)}$, we see that $U(a/f(a)F) = 1$. We deduce  from this that $U(s) = 1$. As $s$ is $\calC$-internal, this implies that $\mathrm{RM}(s) = 1$, and as $F$ is algebraically closed, the type $s$ is strongly minimal. 
    
    Fact \ref{fact: action-classification} rules out the cases of the additive group and an elliptic curve for $\bg{F}{s}{\calC}$. By Lemma \ref{lem: binding iso quotients}, the groups $\bg{F}{p}{\calC}$ and $\bg{F}{s \otimes f(p)}{\calC}$ are isogenous. The projection map $s \otimes f(p) \rightarrow s$ yields a definable surjective morphism $\bg{F}{s \otimes f(p)}{\calC} \rightarrow \bg{F}{s}{\calC}$. As the former is isogenous to a semi-abelian variety, it cannot have any definable morphism to $G_a \ltimes G_m$ nor $\mathrm{PSL_2}$, ruling out the other two cases. Thus $\bg{F}{s}{\calC}$ is definably isomorphic to $G_m(\calC)$, and by connectedness, we get $\bg{f(a)F}{\tp(b/f(a)F)}{\calC} = \bg{F}{s}{\calC}$.

    So $\bg{F}{f(p)}{\calC}$ is definably isomorphic to $E(\calC)$ and $\bg{F}{s}{\calC}$ is definably isomorphic to $G_m(\calC)$.   Corollary \ref{cor: wo from groups} implies that $\bg{F}{s}{\calC, f(p)} = \bg{F}{s}{\calC}$, from which we deduce:    
    \begin{align*}
        \bg{F}{f(p) \otimes s}{\calC , f(p)} & = \bg{F}{s}{\calC, f(p)} \\
        & = \bg{F}{s}{\calC}
    \end{align*}   
    As $f(p)$ is weakly orthogonal to $\calC$, Corollary \ref{cor: wo from groups} also implies that $f(p)$ is weakly orthogonal to $\{ s, \calC \}$. Therefore the short exact sequence:
    \[1 \rightarrow \bg{F}{s}{\calC} \rightarrow \bg{F}{f(p)\otimes s}{\calC} \rightarrow \bg{F}{f(p)}{\calC} \rightarrow 1\]
    \noindent is $F$-definably split by Lemmas \ref{lem: split-iff-retract} and \ref{lem: fund-direct-split-retract}.

    By Lemma \ref{lem: binding iso quotients}, the groups $\bg{F}{s \otimes f(p)}{\calC}$ and $\bg{F}{p}{\calC}$ are isogenous, let $G$ be the $F$-definable group witnessing it. The kernels of the maps to $G$ are given by elements fixing $p$ (resp. $f(p) \otimes s$), and therefore must in particular be contained in the kernels of the maps to $\bg{F}{f(p)}{\calC}$, i.e. $L$ and $\bg{F}{s}{\calC}$. Therefore the two groups are also isogenous, and we have an $F$-definable group $H$ and a commutative diagram of $F$-definable maps:
    \begin{center}
    \begin{tikzcd}
       1 \arrow[r] &  L \arrow[r] \arrow[d] & \bg{F}{p}{\calC} \arrow[r, "\wt{f}"] \arrow[d] & \bg{F}{f(p)}{\calC} \arrow[r] \arrow[d, equal] & 1 \\
       1 \arrow[r] & H \arrow[r] & G  \arrow[r] &  \bg{F}{f(p)}{\calC} \arrow[r] & 1 \\
       1 \arrow[r] & \bg{F}{s}{\calC} \arrow[r] \arrow[u]&  \bg{F}{ s \otimes f(p)}{\calC} \arrow[r] \arrow[u] & \bg{F}{f(p)}{\calC} \arrow[r] \arrow[u, equal] & 1 \\ 
    \end{tikzcd}
    \end{center}
    Both $\bg{F}{p}{\calC}$ and $\bg{F}{f(p) \otimes s}{\calC}$ are abelian, and therefore $F$-definably isomorphic to the $\calC$-points of algebraic groups (see \cite[Lemma 2.2]{jaoui2022abelian} for a proof of that well-known fact). Denote the $\calC$-points of these algebraic groups by $\overline{\bg{F}{p}{\calC}}$ and $\overline{\bg{F}{f(p) \otimes s}{\calC}}$. They are defined over $F \cap \calC$, are isogenous, and by the same reasoning as in the previous paragraph, we obtain an $F \cap \calC$-definable algebraic group $\overline{G}$ and a commutative diagram:

    \begin{center}
    \begin{tikzcd}
       1 \arrow[r] &  G_m(\calC) \arrow[r] \arrow[d, "\iota"] & \overline{\bg{F}{p}{\calC}} \arrow[r] \arrow[d] & E(\calC) \arrow[r] \arrow[d, equal] & 1 \\
       1 \arrow[r] & G_m(\calC) \arrow[r] & \overline{G}  \arrow[r] &  E(\calC) \arrow[r] & 1 \\
       1 \arrow[r] & G_m(\calC) \arrow[r] \arrow[u, "\zeta"]&  \overline{\bg{F}{ s \otimes f(p)}{\calC}} \arrow[r] \arrow[u] & E(\calC) \arrow[r] \arrow[u, equal] & 1 \\ 
    \end{tikzcd}
    \end{center}
    
    \noindent where everything is defined over $\calC \cap F$.

    By Fact \ref{fact: ext-is-functor}, the maps $\iota$ and $\zeta$ give rise to two automorphisms $\iota_*$ and $\zeta_* $ of $ \mathrm{Ext}(E,G_m)$. We obtain that $\iota_*(\overline{\bg{F}{p}{\calC}}) = \overline{G} = \zeta_*(\overline{\bg{F}{ s \otimes f(p)}{\calC}})$.

    We know that $\overline{\bg{F}{ s \otimes f(p)}{\calC}}$ is a $\calC \cap F$-definably split extension of $E(\calC)$, therefore it is trivial as an element of $\mathrm{Ext}(E, G_m)$, and so is its image $\overline{G}$, hence $\overline{G}$ is the identity of $\ext (E,G_m)$, and thus the extension $\overline{\bg{F}{p}{\calC}}$ belongs to the kernel of $\iota_*$. By Proposition \ref{prop: ext-map-fin-ker}, the map $\iota_*$ has finite kernel. Therefore the extension $\overline{\bg{F}{p}{\calC}}$, which is isomorphic to the extension $D$, belongs to a finite subgroup of $\ext (E,G_m)$. This contradicts our choice of $D$, as it was assumed that it did not belong to any finite subgroup of $\ext (G_m,E) = \widehat{E}$. 
\end{proof}

To answer Jin and Moosa's question, we must have a pullback by the logarithmic derivative of a strongly minimal type with binding group an elliptic curve. 

\begin{corollary}\label{cor: jin-moosa-answer}
    There exists an algebraically closed field $F < \mathcal{U}$, a $\calC$-internal type $p \in S_{1}(F)$ such that $\dl : p \rightarrow \dl(p)$ does not almost-split, the type $\dl(p)$ is strongly minimal and has a binding  group isomorphic to the constant points of an elliptic curve $E$.
\end{corollary}

\begin{proof}
    Applying Lemma \ref{lem: reduce to dl} to the algebraically closed field $F$ and $p \in S(F)$ obtained in Theorem \ref{theo: non-split-pullback}, we obtain a type $p' \in S_1(F)$, internal to $\calC$, a fibration $\dl : p' \rightarrow \dl(p')$ with $\dl(p') \subset \dcl(f(p)(\calU)F)$. Because $f$ does not almost split, the map $\dl : p' \rightarrow \dl(p')$ does not almost split either. To conclude, we just need to show that $\dl(p')$ is strongly minimal and $\bg{F}{\dl(p')}{\calC}$ definably isomorphic to the constant points of an elliptic curve.

    Because $\dl(p') \subset \dcl(f(p)(\calU)F)$, we have that $\bg{F}{\dl(p')}{f(p),\calC} = \{ \id \}$, therefore by Lemma \ref{lem: binding iso quotients} we obtain an $F$-definable surjective map $\bg{F}{f(p)}{\calC} \rightarrow \bg{F}{\dl(p')}{\calC}$. As $\bg{F}{f(p)}{\calC}$ is definably isomorphic to $E(\calC)$, and thus strongly minimal, its kernel is either the whole group, which would imply that $\dl(p')$ is $\calC$-definable (i.e. $\dl(p') \subset \dcl(\calC, F)$), or finite. 
    
    Note that the type $\dl(p')$ cannot be algebraic, as if it was, then $\dl : p' \rightarrow \dl(p')$ would be almost split and thus so would $f$. Hence in the first possibility, we deduce that $\dl(p')$ is not weakly orthogonal to $\calC$. But any $\dl(b) \models \dl(p')$ is in the algebraic closure of some $a \models p$, which contradicts $p$ being weakly orthogonal to $\calC$.
    
    Therefore the kernel of the map is finite, and $\bg{F}{\dl(p')}{\calC}$ is definably isomorphic to an algebraic group isogenous to $E(\calC)$, which thus must be an elliptic curve. Note that by the proof of Lemma \ref{lem: reduce to dl}, for any $\dl(b) \models p'$, there is $f(a) \models f(p)$ such that $\dl(b) \in \dcl(f(a),F)$. This, and the previous discussion, forces $\dl(p')$ to be strongly minimal.
\end{proof}

It is unfortunate that our method does not yield a specific differential equation with $\calC$-internal set of solutions. Maybe such an equation could be found using the logarithmic derivative of a non-split semi-abelian variety.

\bibliography{biblio}
\bibliographystyle{plain}

\end{document}